\documentclass[preprint,12pt]{elsarticle}

\usepackage{graphicx}
\graphicspath{{figures/}}
\usepackage{amssymb}
\usepackage{lineno}
\usepackage{amsmath}
\usepackage{stmaryrd}
\usepackage{amsthm}
\usepackage{amsmath}
\usepackage{amssymb}
\usepackage{color}
\usepackage{algorithm}% http://ctan.org/pkg/algorithms
\usepackage{algpseudocode}% http://ctan.org/pkg/algorithmicx
\usepackage{url}
\usepackage{mathtools}

\theoremstyle{definition}
\newtheorem{mydefinition}{Definition} 
\newtheorem{myalgorithm}[mydefinition]{Algorithm} % same counter for Definition and Algorithm

\newcommand{\real}{\mathbb{R}}
\newcommand{\frakA}{\mathfrak{A}}
\newcommand{\frakB}{\mathfrak{B}}
\newcommand{\frakI}{\mathfrak{I}}
\newcommand{\frakS}{\mathfrak{s}}
\newcommand{\frakR}{\mathfrak{r}}

\newcommand{\GridCells}{\mathfrak{K}_h}

\newcommand{\charfunc}{{\chi}}
\newcommand{\jump}[1]{\left\llbracket {#1} \right\rrbracket}
\newcommand{\mean}[1]{\left\{\!\left\{{#1}\right\}\!\right\}}

\newcommand{\dV}{\ \mathrm{dV} }

\newcommand{\divergence}[1]{{\mathrm{div}\left({#1}\right)}}

\newcommand{\DG}[1]{{\mathbb{P}_{#1}}}
\newcommand{\XDG}[1]{{\mathbb{P}^{\text{X}}_{#1}}}
\newcommand{\supp}[1]{\mathrm{supp}{ \left( { #1 } \right)}}

\newcommand{\dS}{{ \ \mathrm{dS} }}

\newcommand{\GammaInt}{{\Gamma_{\mathrm{int}}}}
\newcommand{\GammaDiri}{{\Gamma_\mathrm{D}}}

\newcommand{\innen}{{\mathrm{in}}}
\newcommand{\aussn}{{\mathrm{out}}}

\newcommand{\normI}{{\vec{n}_{\frakI}}}

\newcommand{\normGamma}{{\vec{n}_{\Gamma}}}
\newcommand{\normGammaI}{{\vec{n}_{\frakI,\Gamma}}}

\newcommand{\nablah}{{\nabla_{h}}}

\newcommand{\calC}{\mathcal{C}}
\newcommand{\graph}{{\mathrm{\it{Gr}}}}
\newcommand{\edges}{{\mathrm{\it{Edg}}}}
\newcommand{\aggmesh}{{\mathrm{\it{Agg}}}}
\newcommand{\nI}{\vec{n}_\frakI}
\newcommand{\nOmega}{\vec{n}_{\partial \Omega } }
\newcommand{\gvec}[1]{\underline{#1}}
\newcommand{\gmat}[1]{\underline{\underline{#1}}}

\newcommand{\setV}{\mathbb{V}}
\journal{Computers \& Mathematics with Applications }

\begin{document}

\begin{frontmatter}



\title{BoSSS: a package for multigrid extended discontinuous Galerkin methods}


\author{Florian Kummer and Martin Smuda and Jens Weber}

\address{TU Darmstadt, Darmstadt, Germany}

\begin{abstract}
The software package BoSSS serves the discretization of 
(steady-state or time-dependent)
partial differential equations with discontinuous coefficients
and/or time-dependent domains
by means of an eXtended Discontinuous Galerkin (XDG, resp. DG) method, 
aka. cut-cell DG, aka. unfitted DG.
This work consists of two major parts:
First, the XDG method is introduced and a formal notation is developed, which captures important numerical details 
such as cell-agglomeration and a multigrid framework.
In the second part, iterative solvers for extended DG systems are presented and 
their performance is evaluated.
\end{abstract}


\end{frontmatter}

%\linenumbers


% ===============================================================================
% ===============================================================================
\section{Introduction and Motivation}
% ===============================================================================
% ===============================================================================

BoSSS\footnote{
\url{https://github.com/FDYdarmstadt/BoSSS}}
(\emph{Bo}unded \emph{S}upport  \emph{S}pectral \emph{S}olver) is
a flexible framework for the development, evaluation and application of
numerical discretization schemes 
for partial differential equations (PDEs)
based on the eXtended Discontinuous Galerkin (XDG)
method. 
The main focus of the XDG method are problems with dynamic interfaces and/or time-dependent 
domains.

One important field of application for the XDG method are moving geometries:
Almost all numerical methods for PDEs require the labor-intensive and thus expensive (semi-manual) creation of numerical grids or meshes.
This is particularly cumbersome if time-dependent domains are of interest, e.g. flows with moving boundaries, as can be found when dealing with 
a moving geometry (e.g. rotors, pistons but also flexible structures such as heart valves).
As a side-product, one could use XDG in order to embed a static geometry, which 
would be cumbersome to mesh with higher order elements, on a Cartesian background mesh.

A second field of application are multiphase flows, like mixtures of oil and water: due to topology changes,
e.g. merging or breakup of droplets, such scenarios would be very difficult to handle 
with some conformal meshing or mesh-motion techniques like 
arbitrary Lagrangian-Eulerian (ALE) meshes.

\paragraph{Purpose and scope of this work}
All definitions and algorithms found in this paper have their direct counterpart in the source code.
However, this work is not intended to be a software manual -- in such,
the link to the corresponding place in the source code would have to 
be established with a suitable kind of cross-reference.
Since the source code is constantly being developed, 
the validity of such references would be short-lived.

Instead, this work is intended to serve interested readers as an
\emph{axiomatic description of the extended discontinuous Galerkin framework in BoSSS},
in order to assess its capabilities and limitations.
We therefore emphasize on issues which are not discussed in publications yet,
in specific, details on \emph{multigrid solvers and their performance} will be presented.

% ===============================================================================
\subsection{Prototype Problems}
% ===============================================================================

\paragraph{Poisson with a jump in coefficients}
Let $\Omega \subset \real^D$ ($D \in \{ 2, 3 \}$) be some polygonal domain,
which is partitioned into two disjoint but adjacent sub-domains $\frakA$ and $\frakB$.
These might be referred to as \emph{phases}, since they represent e.g. the liquid and
gaseous phase in multiphase flow simulations.
We assume that the \emph{interface}  $\frakI := \overline{\frakA} \cap \overline{\frakB}$ 
 is a $(D-1)$ - dimensional manifold which is at least 
$\mathcal{C}^1$ almost everywhere.
Therefore, one can define an oriented normal field $\nI$ on $\frakI$; w.l.o.g., we assume
that $\nI$ ``points from $\frakA$ to $\frakB$, i.e. $\nI$ is, on $\frakI,$ equal to 
the outer of $\frakA$ and an inner normal of $\frakB$.
Then, for any property which is continuous in $\Omega \setminus \frakI$, 
one can define the \emph{jump operator}
\begin{equation}
\jump{f}(\vec{x}) :=
    \lim_{\alpha \searrow 0} \left( 
         f(\vec{x} + \alpha \nI) -  f(\vec{x} - \alpha \nI) 
    \right) ,
\label{eq:Jump-Operator-Def}
\end{equation}
for $\vec{x} \in \frakI$.
Then, we can define the piece-wise Poisson problem
\begin{equation}
\left\{ \begin{array}{rll}
- \mu \Delta u                   & = f               & \text{ in } \Omega \setminus \frakI , \\
\jump{u}                         & = 0               & \text{ on } \frakI ,                  \\
\jump{\mu \nabla u \cdot \nI}    & = 0               & \text{ on } \frakI ,                  \\
\vec{u}                          & = g_\text{Diri}   & \text{ on } \Gamma_\mathrm{Diri} ,    \\
\nabla u \cdot \nOmega           & = g_\text{Neu}    & \text{ on } \Gamma_\mathrm{Neu} .     \\
\end{array}
\right.
\label{eq:poisson-jump-problem-def}
\end{equation}
with a discontinuous diffusion coefficient 
\begin{equation}
\mu (\vec{x}) = 
\left\{ \begin{array}{ll}
\mu_\frakA & \text{for } \vec{x} \in \frakA, \\
\mu_\frakB & \text{for } \vec{x} \in \frakB. \\
\end{array} \right.
\label{eq:DiscDiffKoeff}
\end{equation}

\paragraph{Incompressible multiphase flows}
In a transient stetting, the phases are usually considered to be time-dependent
too, i.e. one has the decomposition
$\Omega = \frakA(t) \cup \frakI(t) \cup \frakB(t)$.
Using the same notation as introduced above,
we introduce the incompressible two-phase Navier-Stokes equation for material
interfaces as
\begin{equation}
\left\{ \begin{array}{rll}
   \partial_t \rho \vec{u} 
   +
   \divergence{ \rho \vec{u} \otimes \vec{u} }
   +
   \nabla p
   - 
   \divergence{\mu ( \nabla \vec{u} + (   \nabla \vec{u})^T ) }  
   & = - \rho \vec{G}  
& \text{ in } \Omega \setminus \frakI(t) , \\
   \textrm{div}(\vec{u})  &=  0    
& \text{ in } \Omega \setminus \frakI(t) , \\
   \jump{\vec{u}} & =  0
& \text{ on } \frakI(t) \\ 
    \jump{
      p \nI
      - \mu ( \nabla \vec{u} + (   \nabla \vec{u})^T ) \cdot \nI
    }
     & = 
    \sigma \kappa \normI
& \text{ on } \frakI(t), \\ 
    u  & = \vec{u}_\mathrm{Diri}   
& \text{ on } \Gamma_\mathrm{Diri} ,    \\
  p \nI
      - \mu ( \nabla \vec{u} + (   \nabla \vec{u})^T ) \cdot  \nI           & = 0  
& \text{ on } \Gamma_\mathrm{Neu} .     \\
\end{array} \right.
\label{eq:TwoPhaseNSE}
\end{equation}
with piece-wise constant density and viscosity for both phases, i.e.
\begin{equation}
    \rho(\vec{x}) = \left\{ \begin{array}{ll}
       \rho_\frakA & \textrm{for } \vec{x} \in \frakA \\
       \rho_\frakB & \textrm{for } \vec{x} \in \frakB \\
    \end{array} \right.
    \quad \textrm{and} \quad
    \mu(\vec{x}) = \left\{ \begin{array}{ll}
       \mu_\frakA & \textrm{for } \vec{x} \in \frakA \\
       \mu_\frakB & \textrm{for } \vec{x} \in \frakB \\
    \end{array} \right.
 .
\label{eq:defRhoAndMu}
\end{equation}
Furthermore, $\sigma$ denotes surface tension and $\kappa$ denotes the 
mean curvature of $\frakI$. In a two-phase setting, we also assume 
that the interface does not touches the boundary, 
i.e. $\frakI(t) \cap \partial \Omega = \emptyset$.
If this is not the case, a three-phase contact line occurs,
where two fluids and the solid boundary met at a $D-2$ -- dimensional manifold.
This usually requires more sophisticated boundary conditions (see below).

\subparagraph{Interface tracking and Level-Set:}
Note that problem (\ref{eq:TwoPhaseNSE}) is incomplete without a specification
of the interface motion.
In order to track the individual domains, a 
sufficiently smooth Level-Set field $\varphi$
is introduced; for sake of simplicity, $\varphi$ might also be called `Level-Set'.
Then, time-dependent domains $\frakA(t)$, $\frakB(t)$ and the interface $\frakI(t)$ 
can be described as 
\begin{eqnarray}
\frakA(t) :=  & \left\{ \vec{x} \in \Omega; \ \varphi(t,\vec{x}) < 0 \right\} , \\
\frakI(t) :=  & \left\{ \vec{x} \in \Omega; \ \varphi(t,\vec{x}) = 0 \right\} , \\
\frakB(t) :=  & \left\{ \vec{x} \in \Omega; \ \varphi(t,\vec{x}) > 0 \right\} .
\end{eqnarray}
On $\frakI(t)$, $\varphi$ must therefore comply with the Level-Set equation
\begin{equation}
 \partial \varphi_t + \nabla \varphi \cdot \vec{u} = 0 ,
\label{eq:LevelSetEq} 
\end{equation}
which states that the interface speed in normal direction is equal to 
the flow velocity in normal direction.
From the $\varphi$, the normal $\nI$ can be computed as 
\begin{equation}
\nI = \nabla \varphi / | \nabla \varphi |_2 
\label{eq:normal}
\end{equation}
and   the mean curvature can be computed by Bonnets formula as 
\begin{equation}
 \kappa = \divergence{\frac{\nabla \varphi}{| \nabla \varphi  |_2}} .
\label{eq:Bonnet}
\end{equation}
If the Level-Set is prescribed, e.g. in the case of some externally forced motion, one can infer
the interface speed in normal direction from Eq. (\ref{eq:LevelSetEq}) as
\begin{equation}
 s = \frac{ -\partial_t \varphi}{ | \nabla \varphi | } .
\label{eq:LevSetSpeed}
\end{equation}

\paragraph{Other applications}
The XDG framework can be seen as a multi-purpose technology.
One obvious application is the use as an immersed boundary method (IBM),
e.g. in order to circumvent the (labour-extensive) meshing of domains with complex geometrical details
or to represent domains with moving parts.
Boundary conditions are weakly enforced at the interface $\frakI(t)$ in the same formulation as it would be for boundary fitted methods. 
The fluid domain is represented by $\frakA(t)$, whereas $\frakB(t)$ just describes a void region.
This was demonstrated for e.g. for moving body flows where the motion of solid particles is characterized 
by the Newton-Euler equations
in the work of Krause and Kummer \cite{KrauseKummer2017}.

Immersed boundaries can also be used in the context of compressible flows. 
For inviscid flows around immersed boundaries, such as cylinders and airfoils,
 M\"uller et. al. \cite{MuellerEtAl2017} demonstrated the expected high order of convergence 
 for such geometries using an embedded boundary described by a static Level-Set.

A further issue in the scope of  multiphase flows (\ref{eq:TwoPhaseNSE}) is the interaction of the interface with the domain boundary, 
i.e. the case $\partial \Omega \cap \frakI(t) =: \mathfrak{L}(t) \neq \emptyset$:
this is referred to as the three-phase contact line between fluids $\frakA$ and $\frakB$ and the solid.
In order to allow motion of the \emph{contact line} $\mathfrak{L}(t)$, the no-slip boundary condition has to be relaxed,
yielding the so-called Navier-Slip boundary condition.
%\begin{subequations}
%	\begin{align}
%	\vec{u} \cdot \vec{n}_{\Gamma_S} &= 0,\\
%	\mu \mathbf{P}_S (\nabla{\vec{u}} + \nabla{\vec{u} }^T)  \vec{n}_{\Gamma_S} &= - \beta_S \mathbf{P}_S \vec{u}.
%	\end{align}
%	\label{eq:NavierSlipBC}
%\end{subequations}
%where $\mathbf{P}_S := \mathbf{I} - \vec{n}_{\Gamma_S} \otimes \vec{n}_{\Gamma_S}$. 
It is notable that the XDG-implementation does not need further manipulation 
of the contact line velocity $\vec{U}_{\mathfrak{L}} = (\vec{u} \cdot \vec{n}_{\mathfrak{L}}) \vec{n}_{\mathfrak{L}}$ or contact angle $\Theta$. 
Furthermore, XDG allows the implementation of the generalized Navier-slip boundary condition with 
\begin{equation}
\sigma (\cos{\Theta_{\text{stat}}} - \cos{\Theta}) \vec{n}_{\mathfrak{L}} = \beta_{\mathfrak{L}} (\vec{u} \cdot \vec{n}_{\mathfrak{L}}) \vec{n}_{\mathfrak{L}},
\label{eq:GNBC}
\end{equation}
localized at the contact line $\mathfrak{L}(t)$, which is a $D-2$ -- dimensional manifold.

However, for sake of simplicity and moreover, compactness of the presentation, 
in this work most issues of the XDG method will be discussed 
with respect to the Poisson problem (\ref{eq:poisson-jump-problem-def}). 

% ===============================================================================
\subsection{Historical Development and State of the Art}
% ===============================================================================
The origins of discontinuous Galerkin (DG) methods
 can be tracked back to the work of Reed and Hill \cite{ReedHill1973}.
The name `Discontinuous Galerkin' was mainly established 
through the works of Cockburn and Shu \cite{CockburnAndShu1991},
although similar methods, based on broken polynomial approximation, ideas were already established earlier:
The probably most popular discretization 
for the Laplace operator, the so-called symmetric interior penalty (SIP) method was proposed by Arnold \cite{Arnold1982} about nine years earlier. 

The idea of adapting a finite element method (FEM) to allow jumps in parameters
can be traced back to the 1970s,
where Poisson problems with a jump in the diffusion parameter along a smooth
interface were investigated
by Babu\u{s}ka \cite{Babuska1970}
as well as by Barrett and Elliott \cite{BarrettElliot1987}.
These works mainly relied on isoparametric elements, fitted to the location of the interface at which the discontinuity occurs.
Thus, in a transient setting complex motion of the interface is quite difficult to address.

The continuous extended finite element method (XFEM), was presented by Mo{\"e}s et al. \cite{Moes1999}
to simulate cracking phenomena in solid mechanics. Those ideas were extended to time-dependent
problems by Chessa and Belytschko \cite{ChessaBelytschko2004}.
The first XFEM for two-phase flows was
presented by Gro{\ss} and Reusken \cite{GrosReusken2007,GrosReusken2007b}.
With the so-called intrinsic XFEM presented by Fries \cite{Fries2009},
it became possible to represent also kinks in the velocity field.
This work was later extended by Cheng and Fries \cite{ChengFries2012}
and further by Sauerland and Fries \cite{SauerlandFries2013}.
The first extended DG (XDG) method (also called unfitted DG or cut-cell DG)
was presented by Bastian and Engwer \cite{BastianEngwer2009}
in order to model flows in porous media. Later, those approaches were extended to
multiphase flows by Heimann et. al. \cite{HeimannEtAl2013}.

% ===============================================================================
\subsection{The BoSSS code}
% ===============================================================================

The development of BoSSS has been initiated at the Chair of Fluid Dynamics in 2008.
Since 2017, it is publicly available under the \emph{Apache License}.
The very first motivation for BoSSS is to have a suitable code base for the development of XDG methods.
Very early, it was decided to investigate both, compressible as well as incompressible flows.
The highlights of the code package are:
\begin{itemize}
\item
XDG and support for flows with dynamic interfaces, which is the main topic of this paper.

\item
Suitability for High Performance Computing (HPC): All production algorithms in BoSSS are
implemented MPI parallel. Furthermore, BoSSS provides instrumentation output in order to analyze and
optimize parallel scaling as well as node-level performance. A very important feature  is
dynamic load balancing, which allows re-distribution of the computational mesh when local
processor load changes, e.g. due to motion of the fluid interface.

\item
Rapid prototyping capability: New models resp. equations can be added in a notation that is close to
the usual presentation of DG methods in textbooks, with low development effort. Technical issues, e.g.
the handling of the numerical mesh, are provided by the software library.

\item
Sophisticated workflow and data management facilities: In order to organize and analyze e.g. large
parameter studies, database-centric workflow tools were developed.
\end{itemize}

% ===============================================================================
% ===============================================================================
\section{Discontinuous and Extended Discontinuous Galerkin Methods}
% ===============================================================================
% ===============================================================================

%% ===============================================================================
%\subsection{Discontinuous Galerkin}
%% ===============================================================================

For sake of completeness, and introducing the notation, we briefly summarize 
the DG method; The following definitions are fairly standard and can be found in
similar form in many textbooks \cite{di_pietro_mathematical_2011,hesthaven_nodal_2008}.
\begin{mydefinition}[basic notations] 
We define:
\begin{itemize}
%\item 
%the computational domain: $\Omega \subset \real^D$ which must be polygonal and simply connected; 

\item 
a numerical mesh/grid for $\Omega$ is
a set $\GridCells = \left\{ K_1, \ldots, K_J \right\}$.
% with $h$ being the maximum
%diameter of all cells $K_j$. 
The cells are simply connected, cover the whole domain
($\overline{\Omega} = \bigcup_j \overline{K_j}$),
but do not overlap ($\int_{K_j \cap K_l} 1 \dV = 0$ for $l \neq j$).

\item 
the skeleton of the background mesh:  
$\Gamma := \bigcup_j \partial K_j$. 
Furthermore,  internal edges: $\GammaInt := \Gamma \setminus \partial \Omega$;

\item a normal field $\normGamma$ on $\Gamma$. On $\partial \Omega$, it represents an outer normal,
      i.e., on $\partial \Omega$, $\normGamma = \vec{n}_{\partial \Omega}$.
      By $\normGammaI$, we denote a normal field that is equal to $\normGamma$ on $\Gamma$ and
      equal to $\normI$ on $\frakI$;
\end{itemize}
\end{mydefinition}
For sake of completeness, we notate the approximation space of the 
`standard' DG method:
\begin{mydefinition}[DG space] 
The broken polynomial space of total degree $k$
is defined as:
\begin{equation}
  \DG{k}(\GridCells) := \left\{
    f \in L^2(\Omega); \
    \forall \ K \in \GridCells:
    f|_{K} \textrm{ is polynomial} \right.
    \left. \textrm{and } \mathrm{deg}\left(f|_{K}\right) \leq k
  \right\}
  ;
\end{equation}
\end{mydefinition}

% ===============================================================================
\subsection{Extended Discontinuous Galerkin and Level-Set}
% ===============================================================================
In order to yield a well-defined interface $\frakI$, the Level-Set-field
$ \varphi(t,-)$ must be sufficiently smooth. 
Precisely, we assume $\varphi(t,-)$ to be almost-everywhere in $\calC^1( \Omega )$.
For certain steady-state problems with simple interface shapes, one can 
represent $\varphi$ by an explicit formula that can e.g. be inserted into to code directly.
In more complicated cases, or for temporally evolving problems, $\varphi$ 
may be itself a continuous broken polynomial on the background mesh,
i.e. $ \varphi(t,-) \in \calC^0( \Omega ) \cap \DG{k}(\GridCells) $.
In both cases, one can infer an interface normal (Eq. \ref{eq:normal})
as well as the interface speed (Eq \ref{eq:LevSetSpeed})
from $\varphi$.

\begin{mydefinition}[Cut-cell mesh]
Time-dependent cut-cells are given as
\begin{eqnarray}
K_{j,\frakA} (t) := K_j \cap \frakA(t),
\qquad
K_{j,\frakB} (t) := K_j \cap \frakB(t).
\end{eqnarray}
The set of all cut-cells 
form the \emph{time-dependent cut-cell mesh} $\GridCells^X(t)$.
\end{mydefinition}
Note that we refer to the `original cells' $K_j$ as \emph{background cells}, in contrast 
to the time-dependent cut-cells $K_{j,\frakS} (t)$.
By $\frakS$, resp. $\frakS(t)$ a notation for an arbitrary phase is introduced, 
i.e. $\frakS$ can be either $\frakA$ or $\frakB$.

XDG is essentially a DG method on cut-cells, i.e. one can define the XDG space as follows:
\begin{mydefinition}[XDG space] We define:
\begin{equation}
  \XDG{k}(\GridCells, t) := \DG{k}(\GridCells^X(t))
\end{equation}
\end{mydefinition}
Most DG methods are written in terms of jump- and average operators,
as already defined on $\frakI$, see Eq. (\ref{eq:Jump-Operator-Def}).
This notation is extended onto the skeleton of the mesh $\Gamma$:
\begin{mydefinition}[inner and outer value, jump and average operator]
At the mesh skeleton, the inner- resp. outer-value of a field 
$u \in \mathcal{C}^0(\Omega \setminus \GammaInt \setminus \frakI)$ are defined as:
\[
\begin{array}{rcll}
  u^\innen (\vec{x}) & := &  \lim_{{\xi \searrow 0}} u(\vec{x} - \xi \normGammaI)
  &  \textrm{for } \vec{x} \in \Gamma \cup \frakI \\
  u^\aussn (\vec{x}) & := &  \lim_{{\xi \searrow 0}} u(\vec{x} + \xi \normGammaI)
  &  \textrm{for } \vec{x} \in \GammaInt \cup \frakI \\
\end{array}
\]
Then, the  jump and  average value operator are defined as 
\begin{eqnarray}
\jump{u} & := &
 \left\{  \begin{array}{ll} 
        u^\innen - u^\aussn  &  \text{on } \GammaInt       \\
        u^\innen             &  \text{on } \partial \Omega 
 \end{array} \right. ,
\label{eq:defJumpOperator}
\\
\mean{u} & := &
 \left\{  \begin{array}{ll} 
        (u^\innen + u^\aussn)/2  &  \text{on } \GammaInt       \\
        u^\innen                 &  \text{on } \partial \Omega 
 \end{array} \right.  .
\label{eq:defMeanOperator}
\end{eqnarray}
\end{mydefinition}

%\begin{itemize}
%\item the broken gradient $\nablah$: for $u \in \mathcal{C}^1(\Omega \setminus \Gamma \setminus \frakI)$,
%      $\nablah u$ denotes the gradient on the domain $\Omega \setminus \Gamma \setminus \frakI$;
%      in analog fashion, the broken divergence $\divergenceh{\vec{u}}$;
%
%\item the measure $| X |$ of some set $X$ is given as $| X | = \int_{X} 1 \dV$;
%
%\item the species-volume-fraction $\spcFrac{\frakS}{K}$ for some cell $K \in \GridCells$:
%      $\spcFrac{\frakS}{K} := | K \cap \frakS |  / | K |$,
%      for $\frakS \in \{ \frakA, \frakB \}$;
%
%\item the standard-basis vector $\vec{e}_d$, for $d \in \{ 1, 2 \}$:
%      $\vec{e}_1 = (1,0)$, $\vec{e}_2 = (0,1)$;
%\end{itemize}

For the implementation of an XDG-method, accurate 
numerical integration on the cut-cells $K_{j,\frakS}(t)$ is required.
In BoSSS the user can either select 
the  Hierarchical Moment Fitting (HMF) procedure, developed by M\"uller et. al. \cite{Mueller2013},
or, alternatively,
a method proposed by Saye \cite{Saye2015}.
While the HMF supports all types of cells (triangles, quadrilaterals, tetra- and hexahedrons)
the method of Saye is generally faster, but restricted to quadrilaterals and hexahedrons.

\paragraph{Ensuring continuity of $\varphi$, computation of normals and curvature}
Some XDG method for the two-phase Navier Stokes problem (\ref{eq:TwoPhaseNSE}) has to be coupled
with a second method to compute the evolution of the interface. 
Since the respective Level-Set equation (\ref{eq:LevelSetEq}) is of hyperbolic type,
a conventional DG approximation seems a reasonable choice.
However, such a method will typically yield a discontinuous field $\varphi$.
In such cases, $\varphi$ needs to be projected to a continuous space before it can be used to 
setup the XDG method.
Furthermore, the curvature $\kappa$ required in problem (\ref{eq:TwoPhaseNSE}),
see also Eq. (\ref{eq:Bonnet}), needs to be handled with care.
For a detailed presentation of the filtering procedures used in BoSSS, 
we refer to the work of Kummer and Warburton  \cite{KummerWarburton2016}.
For sake of simplicity, in this work we assume $\varphi$ to be sufficiently smooth
with respect to space and time.

% ===============================================================================
\subsection{Variational formulations}
% ===============================================================================
Now, problems (\ref{eq:poisson-jump-problem-def}) and (\ref{eq:TwoPhaseNSE}) 
can be discretized in the XDG space.
In general, we are interested in systems like the Navier-Stokes equation, 
with $D_v$ dependent variables which are not necessarily 
discretized with the same polynomial degree.
We therefore define the degree-vector $\gvec{k} = (k_1, \ldots k_{D_v})$;
and introduce an abbreviation for the function space of test and trial functions, 
\begin{equation}
\setV_{\gvec{k}}^X(t) := \prod_{\gamma=1}^{D_v}  \XDG{k_\gamma}(\GridCells, t).
\label{eq:defXDGvectorSpace}
\end{equation}
Then, the discrete version of some linear problem,
for a fixed time $t$
 formally read as:
find $U \in \setV_{\gvec{k}}^X(t)$ so that
\begin{equation}
   a(U,V) = b(V) \quad \forall V \in \setV_{\gvec{k}}^X(t) .
\label{eq:GenericVarScheme}
\end{equation}

\paragraph{Discrete variational formulation for Poisson Eq. (\ref{eq:poisson-jump-problem-def})}
The variational formulation of the symmetric interior penalty method,
originally proposed by Arnold \cite{Arnold1982},
reads as
\begin{multline}
a_{\text{sip}}(u,v)
:= 
  - \oint_{\frakI \cup \Gamma \setminus \Gamma_{\text{Neu}}} 
        \mean{\mu \nabla_h u} \cdot \vec{n}_{\frakI,\Gamma} \jump{v}
      + 
        \mean{\mu \nabla_h v} \cdot \vec{n}_{\frakI,\Gamma} \jump{u}
      \dS \\
  + \oint_{\frakI \cup \Gamma \setminus \Gamma_{\text{Neu}}}
       \eta \max\{\mu^\innen, \mu^\aussn \} \jump{u}\jump{v}
     \dS
  + \int_{\Omega} \mu \nabla_h u \cdot \nabla_h v \dV
\label{eq:DefSIP}
\end{multline}
for the left-hand-side of Eq. (\ref{eq:GenericVarScheme}).
Here,  $\nablah u$ denotes the broken gradient, 
where differentiation at the jumps on $ \Gamma \cup \frakI$ is excluded.
The linear form on the right-hand-side of Eq. (\ref{eq:GenericVarScheme}) is given as
\begin{equation}
 b(v) := 
  \int_\Omega f v \dV
 - \oint_{\Gamma_{\text{Diri}}}
              \mu g_{\text{Diri}} \left( \nabla_h v \cdot \vec{n}_{\partial \Omega}  
            - \eta v \right)
      \dS
   +  \oint_{\Gamma_{\text{Neu}}}
              \mu g_{\text{Neu}} v
      \dS
\label{eq:SIPrhs}
\end{equation}
The SIP factor $\eta$ is known to scale as $\eta \approx \text{const} \cdot k^2 / h'$,
where $h'$ is a local length scale of the agglomerated cut-cell, see section \ref{sec:Agglom}.
For certain specific cell shapes, explicit formulas for $\eta$ can be given,
a comprehensive overview is given in the thesis of Hillewaert \cite{HillewaertThesis2013}.
In the case of XDG methods whith arbitrary cell shapes, some rules-of-thumb can be used, 
given that a sufficiently large multiplicative constant is used.
Alternatively, a symmetric weighted interior penalty (SWIP) form,
as presented in the textbook of Di Pietro \& Ern \cite{di_pietro_mathematical_2011}
might be used.

\paragraph{Discrete variational formulation for Navier-Stokes Eq. (\ref{eq:TwoPhaseNSE})}
Due to the non-linearity in the convectional term, the Navier-Stokes system is usually 
split up into the linear Stokes part $a(-,-)$ and 
the nonlinear convection operator $c(-,-,-)$. 
It is known that, in order to obtain an inf-sup stable discretization, 
the DG polynomial degree of the pressure has to be one lower than for velocity,
i.e. velocity is discretized with degree $k$ and pressure with degree $k-1$.
()Note this is proven only for special cases, see the work of Girault et. al. \cite{Girault2005}.)
In this notational framework, for spatial dimension $D=3$, we write $D_v = D + 1$ 
and $\gvec{k} = (k,\ldots,k,k-1)$.
The variational formulation of the Navier-Stokes equation formally reads as:
find 
$
U =( \vec{u}, p) \in \setV_{\gvec{k}}^X(t) 
   = \left( \XDG{k}(\GridCells, t) \right)^D \times \XDG{k-1}(\GridCells, t)
$
so that
\begin{equation}
  \int_\Omega \partial_t U \cdot V \dV
  +
  c(U,U,V) + a(U,V) = b(V) \quad \forall V \in \setV_{\gvec{k}}^X(t) .
\label{eq:NSEvariational}
\end{equation}
A complete specification of the involved forms would be too lengthy here; hence,
we refer to the works of Heimann et. al. \cite{HeimannEtAl2013} and Kummer \cite{KummerXnse2017}.
Obviously, (i) this equation still requires a temporal discretization 
and (ii) a nonlinear solver.

For sake of simplicity, we assume an implicit Euler discretization in time,
i.e. $\partial_t U \approx (U^1 - U^0)/\varDelta t$,
where $U^0 \in \setV_{\gvec{k}}^X(t^0)$ denotes the known value from previous timestep $t^0$,
and $U^1 \in \setV_{\gvec{k}}^X(t^1)$ denotes the unknown value at the new timestep $t^1$,
i.e. $t^1 - t^0 = \varDelta t$.
We also fix the first argument of $c(-,-,-)$. Then, scheme (\ref{eq:NSEvariational})
reduces to a form which is formally equivalent to scheme (\ref{eq:GenericVarScheme}),
namely: find $U^1 \in \setV_{\gvec{k}}^X(t^1)$ so that
\begin{multline}
  \int_\Omega U^1 \cdot V \dV
  +
  c(U_0,U^1,V) + a(U^1,V) =  \\
  b(V) + \int_\Omega U^0 \cdot \mathcal{L}_{t^1}^{t^0} V \dV 
  \quad 
  \forall V \in \setV_{\gvec{k}}^X(t^1) .
\label{eq:NSEvariationalLinearized}
\end{multline}
The linearization point $U_0$ may be either set as $U_0 = U^0$, which results 
in a semi-implicit formulation.
Alternatively, one can utilize fully implicit approach and
iterate over equation (\ref{eq:NSEvariationalLinearized}), so that 
$U_0 \rightarrow U^1$ or employ a Newton method.

On the right-hand-side, the linear operator $\mathcal{L}_{t^1}^{t^0}$ performs a lifting 
from $\setV_{\gvec{k}}^X(t^1)$ to $\setV_{\gvec{k}}^X(t^0)$. 
This requires that the corresponding part of the right-hand-side of (\ref{eq:NSEvariationalLinearized})
must be integrated on the old cut-cell mesh at time $t^0$.

The lifting is defined as follows:
\begin{mydefinition}[temporal XDG space lifting]
\label{def:liftingOp}
The cut-cell mesh $\GridCells^X(t^0)$ and $\GridCells^X(t^1)$ at times $t^0$ and $t^1$ 
have equal topology, if, and only if for each cut-cell one has 
\begin{eqnarray*}
| K_{j,\frakS} (t^0) | > 0  & \quad \Rightarrow \quad & | K_{j,\frakS} (t^0) \cap K_{j,\frakS} (t^1) | > 0  \text{ and } \\ 
| K_{j,\frakS} (t^1) | > 0  & \quad \Rightarrow \quad & | K_{j,\frakS} (t^1) \cap K_{j,\frakS} (t^0) | > 0 . 
\end{eqnarray*}
On meshes with equal topology, the lifting operator  
\[
\XDG{k}(\GridCells, t^1) \ni 
u^1   \mapsto \mathcal{L}_{t^1}^{t^0} u^1 =: u^0 
\in  \XDG{k}(\GridCells, t^0)
\]
is uniquely defined 
requiring polynomial equality on the common domain of old and new cut-cell, 
i.e. by the property 
\[
 \left. u^0 \right|_{ K_{j,\frakS} (t^0) \cap K_{j,\frakS} (t^1)} 
 = 
 \left. u^1 \right|_{ K_{j,\frakS} (t^0) \cap K_{j,\frakS} (t^1)}
\quad 
\forall j, \frakS 
.
\]
\end{mydefinition}
Since $\setV_{\gvec{k}}^X(t)$ is a product of spaces $\XDG{k}(\GridCells, t)$,
the lifting $\mathcal{L}_{t^1}^{t^0} u^1$ naturally extends to $\setV_{\gvec{k}}^X(t)$.
Obviously, it cannot always be ensured that 
the cut-cell meshes for $t^0$ and $t^1$ have equal topology.
However, this important property can be achieved through 
cell agglomeration, which is addressed in the next section.
% ===============================================================================
\subsection{Cell agglomeration}
\label{sec:Agglom}
% ===============================================================================
The motivation for aggregation/agglomeration meshes is three-fold:
removal of small cut-cells,
avoiding topology changes in the cut-cell mesh for a single time-step 
 and formulation of multigrid methods
without the usual hierarchy of meshes, cf. Section \ref{sec:MGsolvers} and Figure \ref{fig:agg_multigrid}).

Formally, the aggregation mesh is introduced by means of graph theory:
\begin{mydefinition}[Graph of a numerical mesh] 
Let $\mathfrak{K}$ be a numerical mesh;
For $K_1, K_2 \in \mathfrak{K}$, the set $\{K_1, K_2\}$ is called a 
\emph{logical edge} if, and only if 
$\oint_{\overline{K_1} \cap \overline{K_2}} 1 \dS > 0$.
Furthermore, let $\edges( \mathfrak{K} )$ be the set of all logical edges.
Then, the pair $( \mathfrak{K}, \edges( \mathfrak{K} ) ) =: \graph(\mathfrak{K})$
forms an undirected graph in the usual sense.
\end{mydefinition}

\begin{mydefinition}[Aggregation maps and meshes]
Let $A \subset \edges( \mathfrak{K} )$ be an aggregation map and
$ a := \{ K_1, \ldots , K_L  \} $ be the nodes of a connected component of  $\graph(\mathfrak{K})$.
Note that $a$ might consist of only a single element, i.e. an isolated node, 
which is called an non-aggregated cell with respect to $A$.
The aggregate cell is defined as the union of all cells, i.e. 
$K_a := \bigcup_{K_l \in a}^\partial K_l $
(Rem.: in order to ensure that the aggregate cell is again a simply connected, open set,
one has to take the closure of each cell first and then subtract the boundary,
therefore we define a modified union as 
$X \cup^\partial Y := 
   ( \overline{X} \cup \overline{Y} ) 
   \setminus \partial 
( \overline{X} \cup \overline{Y} )$.)

For  $A \subset \edges( \mathfrak{K} )$ 
the aggregation mesh
$\aggmesh(\GridCells, A)$ is the set of all aggregate cells which can be formed w.r.t. $A$.
\end{mydefinition}
Based upon the aggregation mesh, an aggregated XDG space can be defined:
\begin{mydefinition}[Aggregated XDG space] 
For some agglomeration map $A \in \edges(\GridCells^X(t^1))$
we define the agglomerated XDG space as:
\begin{equation}
  \XDG{k}^A(\GridCells, t) := \DG{k}(\aggmesh(\GridCells^X(t), A))
\end{equation}
\end{mydefinition}
Obviously, the agglomerated XDG space is a sub-space of the original space, i.e.
$\XDG{k}^A(\GridCells, t) \leq_\real \XDG{k}(\GridCells, t)$. 

\paragraph{Temporal discretization and stabilization against small cut-cells}
As already noted above, in order to discretize temporally evolving systems such 
as Eq. (\ref{eq:NSEvariationalLinearized}), one has to ensure that the cut-cell mesh 
at time steps $t^0$ and $t^1$ have equal topology in order to obtain a well-defined method.
Otherwise, the required lifting operator (see Definition \ref{def:liftingOp}) is undefined.
For multi-step schemes, which involve 
multiple time steps, the topology has to be equal for all time steps; the same holds for the 
intermediate steps of Runge-Kutta schemes, cf. \cite{Kummer2017}.

Furthermore, since the interface position is arbitrary, cut-cells can be arbitrarily small,
i.e. its volume fraction $ | K_{j,\frakS}(t) | / | K_{j} | $ w.r.t. the background cell can be small.
This leads e.g. to large penalty parameters $\eta$ in the SIP form (\ref{eq:DefSIP}) which 
is known to cause undesirably high condition numbers of the discretized system.

Therefore, instead of solving the variational system on the space $\setV_{\gvec{k}}^X(t)$,
which is induced by the cut-cell mesh, one employs an XDG space on an appropriately agglomerated cut-cell mesh.
A valid agglomeration map 
$A_{\alpha,t^1,t^0} \subset \edges(\GridCells^X(t^1)) \cup \edges(\GridCells^X(t^0))$
for these purpose must meet the following requirements:
\begin{itemize}
\item
The meshes $\aggmesh(\GridCells^X(t^1), A_{\alpha,t^1,t^0} )$ and 
           $\aggmesh(\GridCells^X(t^0), A_{\alpha,t^1,t^0} )$
have the same topology.

\item 
All cut-cells with a volume fraction 
$0 < | K_{j,\frakS}(t) | / | K_{j} | \leq \alpha$ are agglomerated

\item 
There is no agglomeration across species, i.e. there exists no edge 
$\{ K_{j,\frakA}, K_{l,\frakB} \}$ in $ A_{\alpha,t^1,t^0} $.
\end{itemize}
The formulation of an algorithm that constructs an agglomeration map which fulfills the
properties noted above is left to the reader. 
It may consist of a loop over all cut-cells. The cut-cell $K_{j,\frakS}$ must be 
agglomerated if it is a new cell (i.e. $| K_{j,\frakS}(t^1) | > 0$ and $| K_{j,\frakS}(t^0) | = 0$)
or a vanished cell (i.e. $| K_{j,\frakS}(t^1) | = 0$ and $| K_{j,\frakS}(t^0) | > 0$) 
or if its volume fraction is below the threshold $\alpha$. A decent value for $\alpha$
lies in the range of 0.1 to 0.3, cf. \cite{KummerXnse2017,MuellerEtAl2017}.
In our implementation, such a cell is agglomerated to its largest neighbor cell in 
the same species.

\paragraph{The final system}
Instead of solving the generic variational system (\ref{eq:GenericVarScheme})
on the space $\setV_{\gvec{k}}^X(t^1)$, the aggregated space
\begin{equation}
\setV_{\gvec{k}}^{X,\alpha,\varDelta t} 
:= 
\prod_{\gamma=1}^{D_v}  \XDG{k}^{ A_{\alpha,t^0,t^1} }(\GridCells, t^1)
\end{equation}
is used for discretization. 
In comparison to (\ref{eq:defXDGvectorSpace}), the dependence on time  $t$
is dropped since w.l.o.g. all temporally evolving systems are solved for 
the `new' time step $t^1$.

Hence, the final discretization reads as:
find $U \in \setV_{\gvec{k}}^{X,\alpha,\varDelta t} $ so that
\begin{equation}
   a(U,V) = b(V) \quad \forall V \in \setV_{\gvec{k}}^{X,\alpha,\varDelta t}  .
\label{eq:GenericVarSchemeAgglommSpace}
\end{equation}

Over the course of multiple time-steps, the agglomeration graph typically changes.
After each time step is complete, the solution is injected into the non-agglomerated space.
For the next time step it is projected back (in an $L^2$-sense) onto the 
(potentially different) agglomerated space to serve as an initial value.

% ===============================================================================
% ===============================================================================
\section{Multigrid solvers}
\label{sec:MGsolvers}
\begin{figure}
\[
\begin{array}{rccc}
&
\textrm{level 1}: \GridCells^1 = \GridCells &
\textrm{level 2}: \GridCells^2  &
\textrm{level 3}: \GridCells^3  \\
\rotatebox{90}{\rule{0.5cm}{0.0cm} cartesian $8 \times 8$} &
\includegraphics[width=0.3\textwidth]{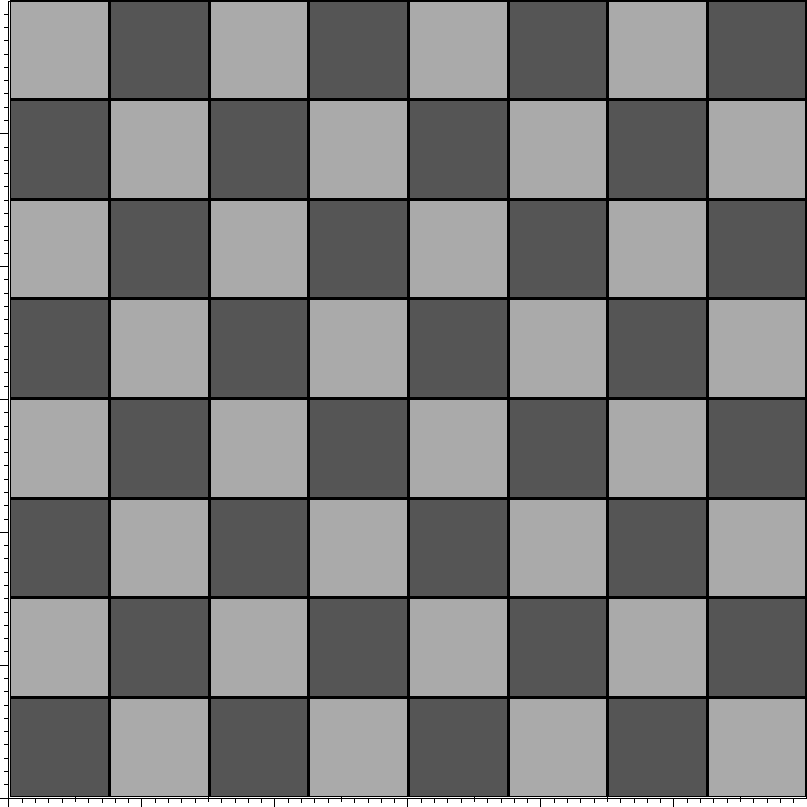} &
\includegraphics[width=0.3\textwidth]{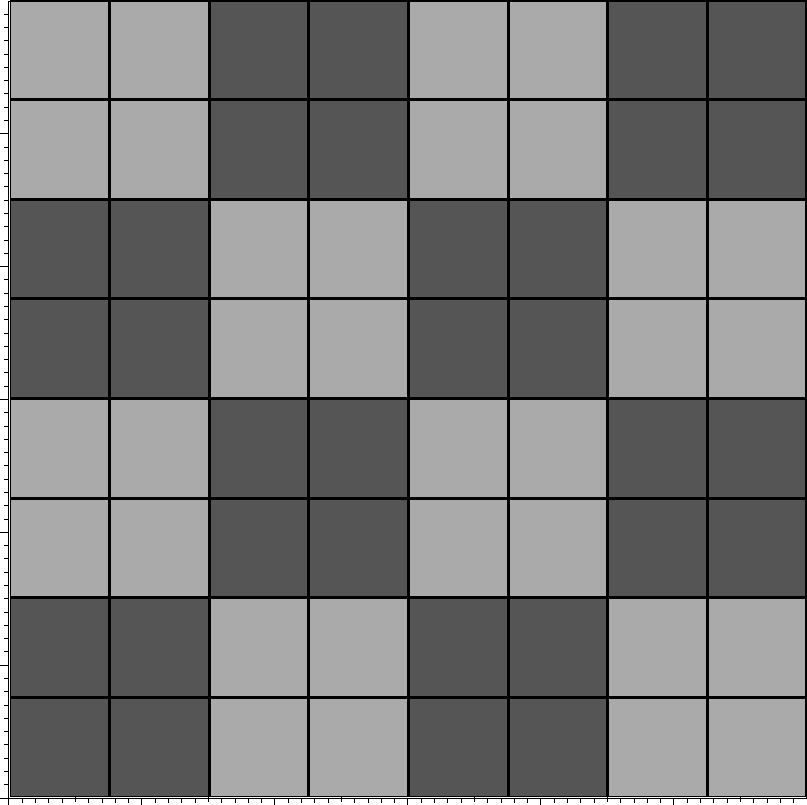} &
\includegraphics[width=0.3\textwidth]{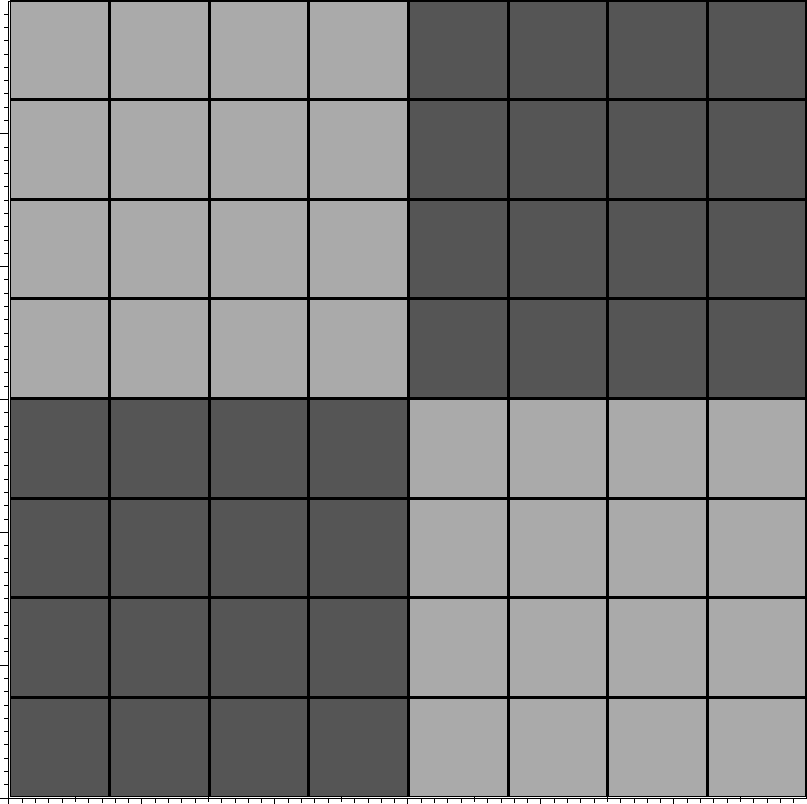} \\
\rotatebox{90}{\rule{0.5cm}{0.0cm} cartesian $9 \times 9$} &
\includegraphics[width=0.3\textwidth]{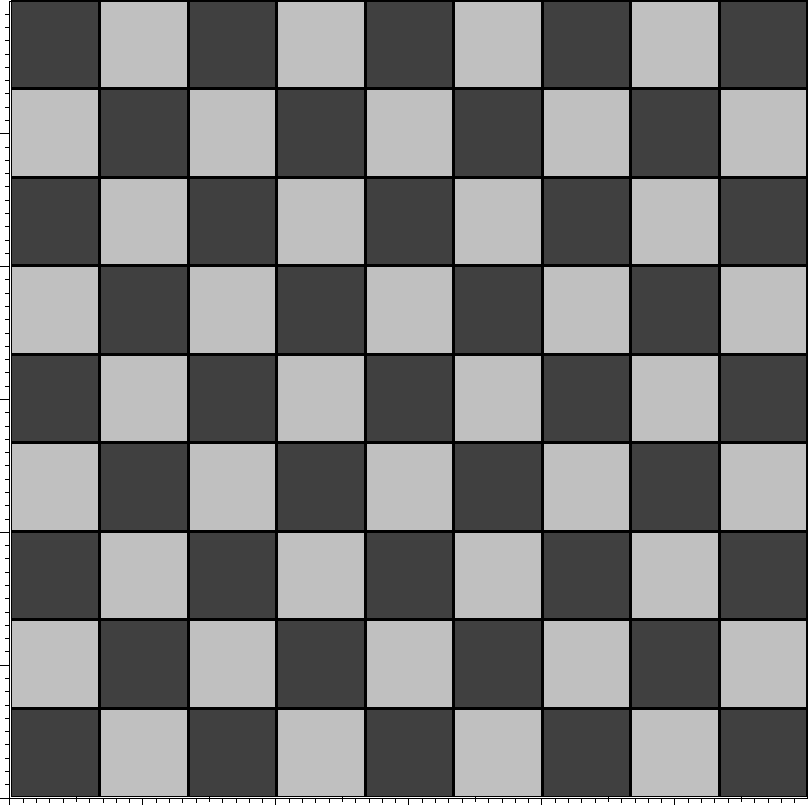} &
\includegraphics[width=0.3\textwidth]{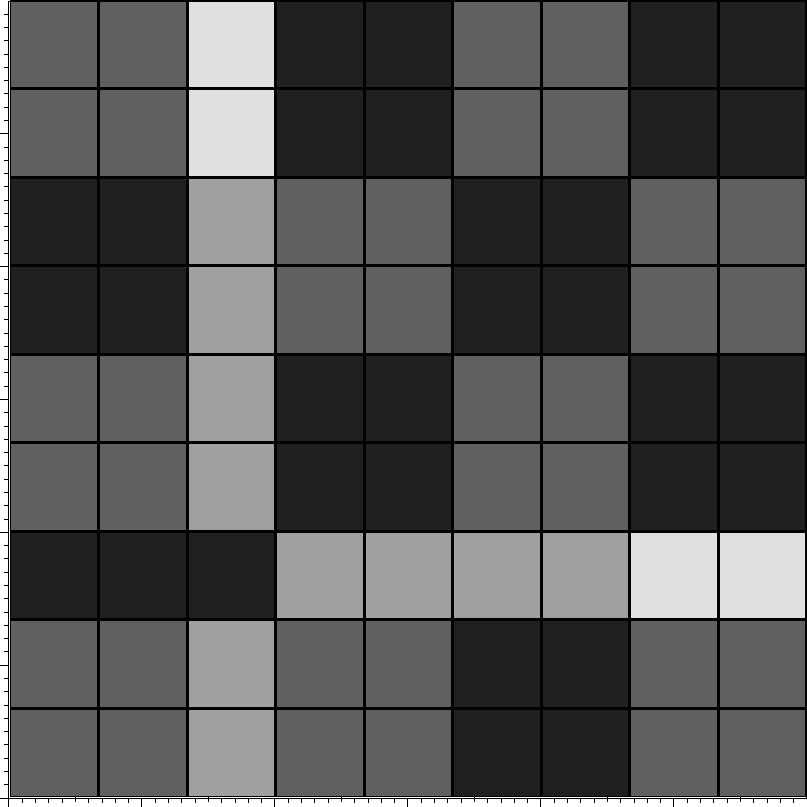} &
\includegraphics[width=0.3\textwidth]{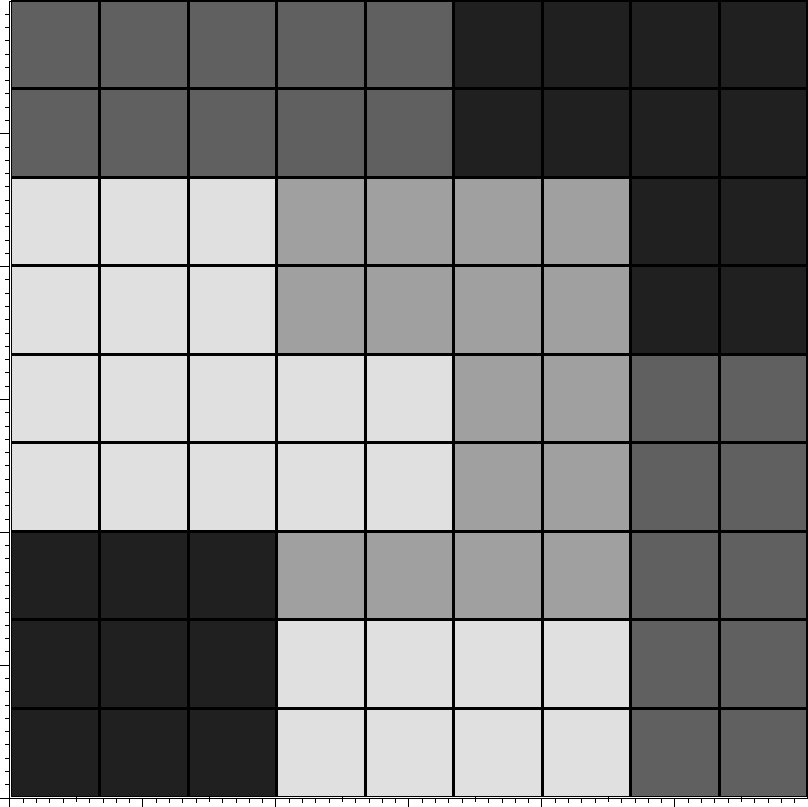} \\
\rotatebox{90}{\rule{0.25cm}{0.0cm} cartesian 'L' $7 \times 7$} &
\includegraphics[width=0.3\textwidth]{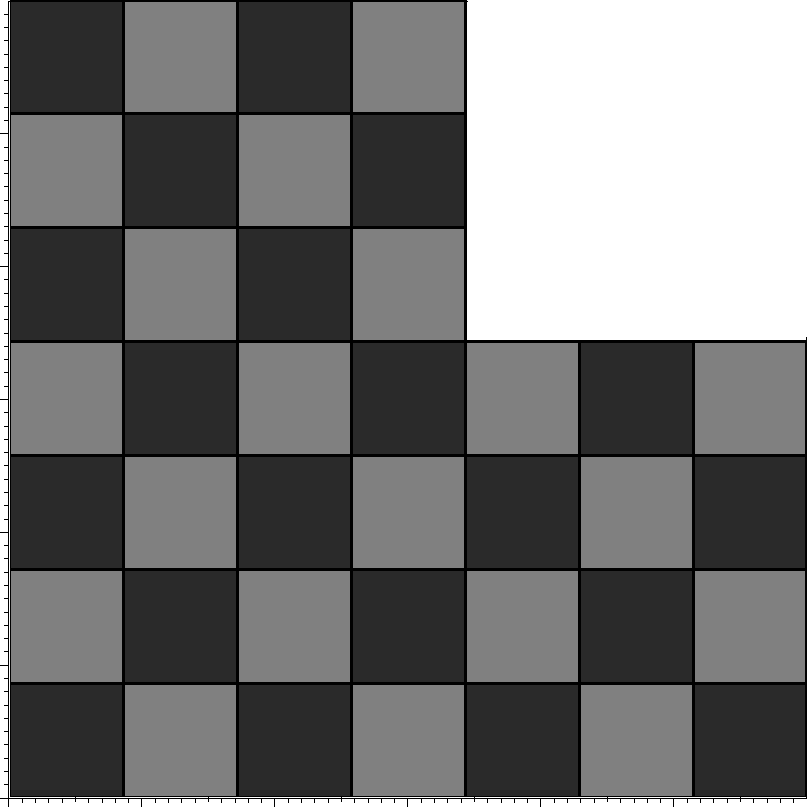} &
\includegraphics[width=0.3\textwidth]{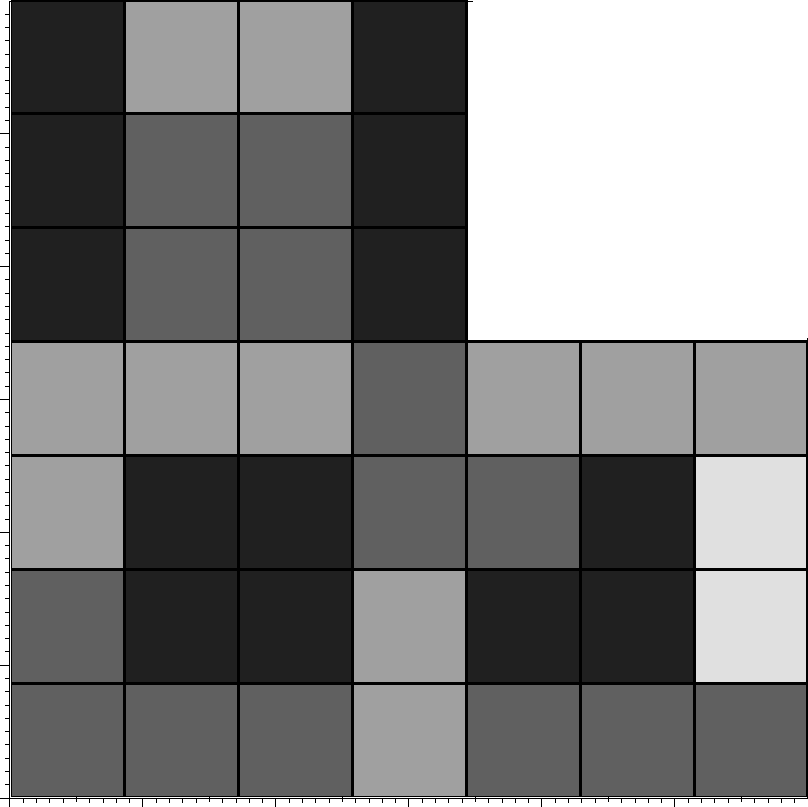} &
\includegraphics[width=0.3\textwidth]{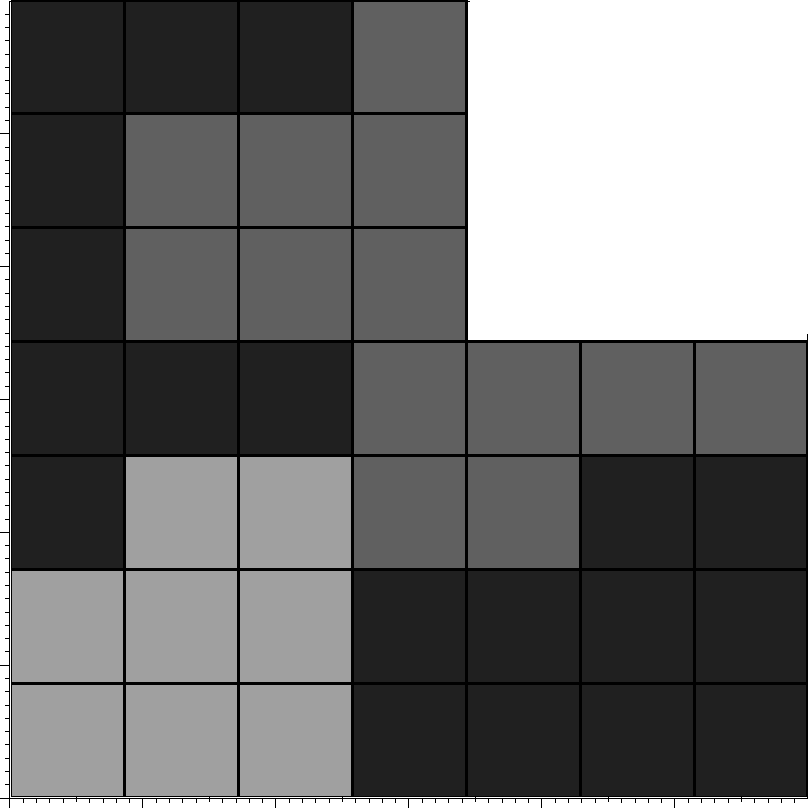} \\
\rotatebox{90}{\rule{0.25cm}{0.0cm} cartesian 'L' $9 \times 9$} &
\includegraphics[width=0.3\textwidth]{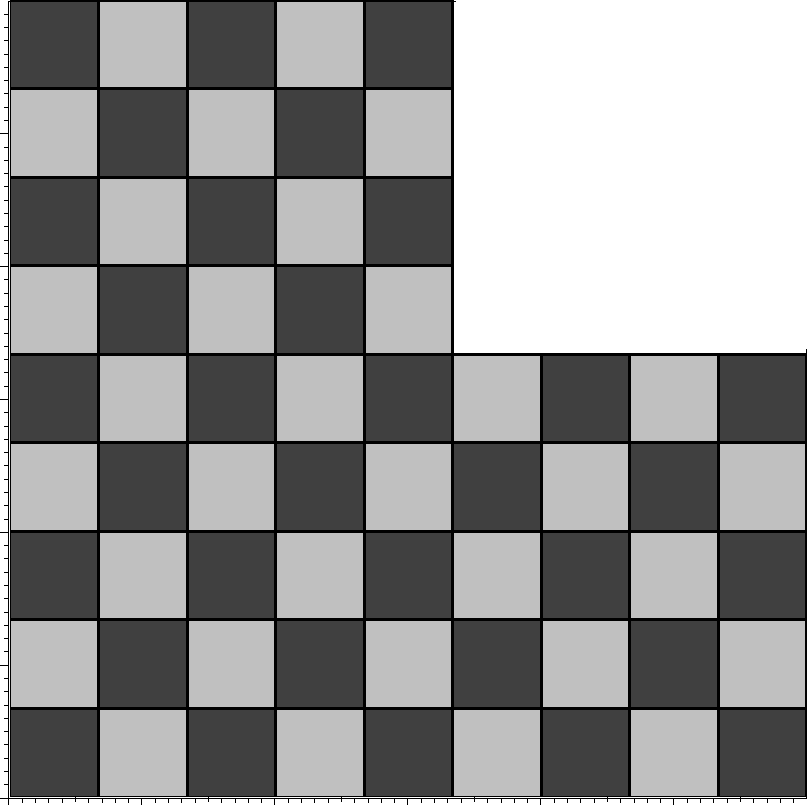} &
\includegraphics[width=0.3\textwidth]{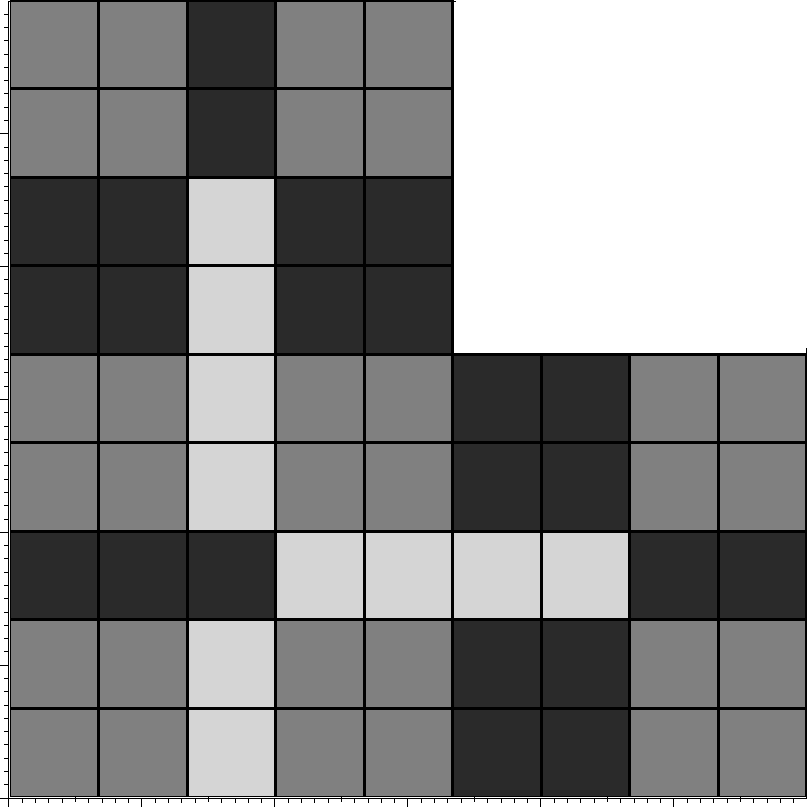} &
\includegraphics[width=0.3\textwidth]{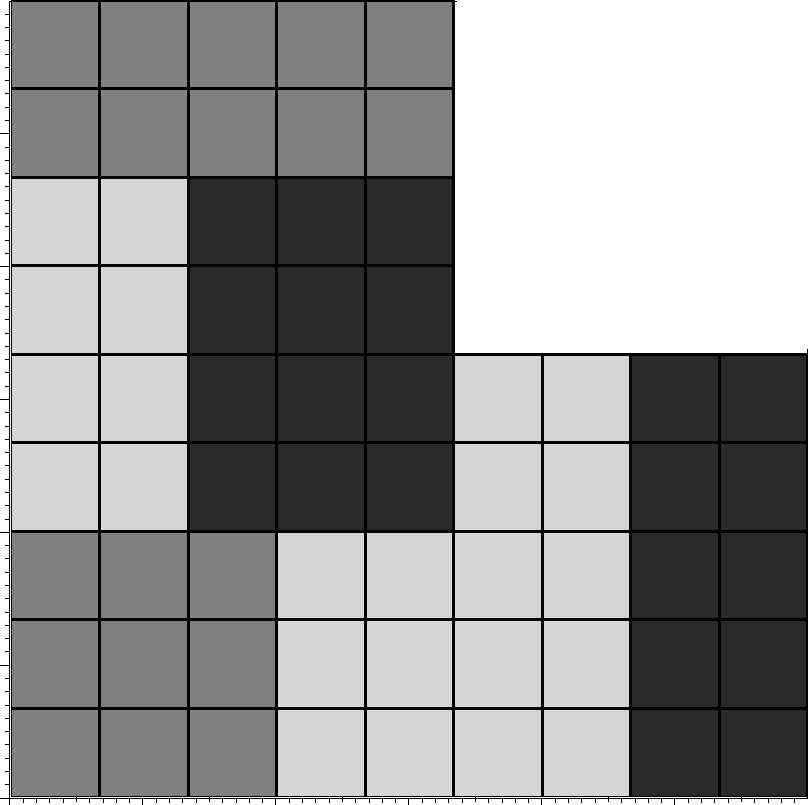} \\
\end{array}
\]
\caption{Examples of aggregation multigrid:
an advantage of this approach is that it is with respect to geometry more flexible
than the classical geometric multigrid approach where the finest grid level
is induced by the coarsest.
}
\label{fig:agg_multigrid}
\end{figure}

% ===============================================================================
% ===============================================================================
The remaining of this paper is dedicated to the solution of the linear system (\ref{eq:GenericVarSchemeAgglommSpace}).
The solvers presented are use a combination of aggregation- and p-multigrid.
Aggregation multigrid can be seen as a combination of conventional h-multigrid
and algebraic multigrid methods: there is still an underlying mesh of polyhedral cells.
Due to these cells, the flexibility is comparable to algebraic multigrid, 
see Figure \ref{fig:agg_multigrid}.

This section is organized as follows: 
first, the aggregation multigrid framework will be unified with the XDG method 
established so far (section \ref{sec:AggMGforXDG}).
Next, the construction of a basis for the nested sequence of approximation
spaces will be laid out, in order tor transfer the basis-free 
formulation  (\ref{eq:GenericVarSchemeAgglommSpace}) into a matrix form.
Finally, the specific combination of multigrid algorithms are presented (\ref{sec:Precond}).

% ===============================================================================
\subsection{Aggregation multigrid for XDG}
\label{sec:AggMGforXDG}
\begin{figure}
%\begin{center}
%\textcolor{green}{Todo: Smuda: nice version of this image}
%\includegraphics[width=0.5\textwidth]{CutCellAggMap}
%\end{center}
\centering
\def\svgwidth{\textwidth}
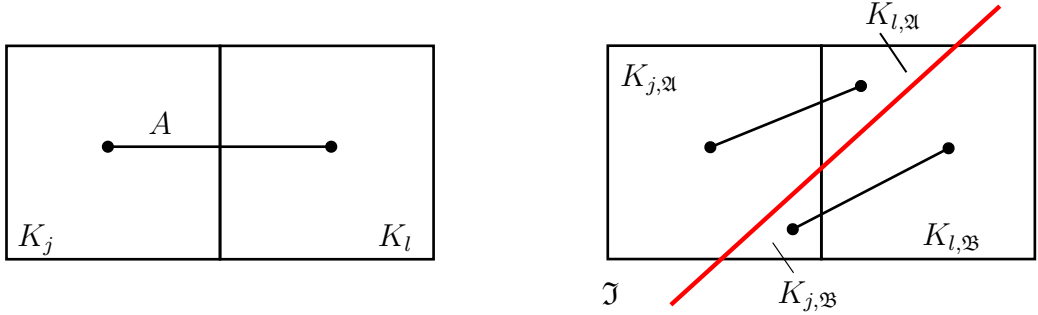
\caption{
From an aggregation map $A \subset \edges(\GridCells)$ (left) to 
the corresponding map  $A^X(t) \subset \edges(\GridCells^X(t))$ on the cut-cell mesh (right).
Since the edges of $A$ are duplicated for each species, no 
aggregation across the interface $\frakI$ occurs.
}
\label{fig:CutCellAggMap}
\end{figure}

% ===============================================================================
The starting point of the aggregation multigrid is a 
sequence of aggregation maps
\begin{equation}
 \emptyset = A^1 \subset A^2 \subset \ldots \subset A^\Lambda \subset \edges( \GridCells )
\label{eq:MultigridAggMeshSeq}
\end{equation}
on the background mesh.
Note that the injection/projection operator between aggregation grid are quite 
expensive to compute. Therefore, it is beneficial to compute them initially and only update 
when necessary in cut-cells. 
Since these are defined on the background mesh, in order to be precomputed,
one cannot directly apply these aggregation maps onto the cut-cell mesh.
Therefore, aggregation maps from the background mesh must be mapped onto the cut-cell mesh:
\begin{mydefinition}[Mapping of an aggregation map onto a cut-cell mesh]
Given is an aggregation map $A \subset \edges(\GridCells)$ on the background mesh $\GridCells$.
The corresponding aggregation map $A^X(t)$ on the cut-cell mesh $\GridCells^X(t)$
is formed from all edges $\{ K_j, K_l \} \in A$ by duplicating them for each species,
i.e. from edges 
\[
  \{ K_{j,\frakS} (t), K_{l,\frakS} (t) \}
   \quad
   \forall \frakS \in \{\frakA, \frakB \}, 
   \text{ if } | K_{j,\frakS} (t) | > 0 \text{ and } | K_{l,\frakS} (t) | > 0 .
\]
\end{mydefinition}
Note that this construction avoids aggregation across the interface $\frakI$,
as illustrated in Figure \ref{fig:CutCellAggMap}.
Through the mapping of an aggregation map $A \mapsto A^X(t)$, 
a sequence of aggregated XDG spaces is induced,
\begin{equation*}
\setV_{\gvec{k}}^{X,\alpha,\varDelta t} =: \setV_{\gvec{k}}^{X,1} 
  \ \prescript{}{\real}\geq \
\setV_{\gvec{k}}^{X,2} 
  \ \prescript{}{\real}\geq \
\ldots
  \ \prescript{}{\real}\geq \
\setV_{\gvec{k}}^{X,\Lambda} 
\end{equation*}
where the space on mesh level $\lambda$ can be defined as 
\begin{equation}
\setV_{\gvec{k}}^{X,\lambda} := 
\prod_{\gamma=1}^{D_v}  \XDG{k}^{ A'^{\lambda} }(\GridCells, t^1) 
\end{equation}
and the aggregation map 
is the union of the predefined multigrid aggregation sequence (see \ref{eq:MultigridAggMeshSeq})
and the aggregation map $A_{\alpha,t^1,t^0}$ to stabilize small cut-cells and prevent temporal toplology changes,
i.e.
\begin{equation}
 {A'}^{\lambda} := (A^\lambda)^X (t^1) \cup A_{\alpha,t^1,t^0} .
\end{equation}
% ===============================================================================
\subsection{Basis representation and indexing}
\label{seg:MGBasis}
% ===============================================================================
Up to this point, the XDG method is formulated in a variational, coordinate-free 
form. In order to notate solver algorithms in a typical 
matrix-notation, a basis representation of the XDG space is required.
\paragraph{A Basis of $\setV_{\gvec{k}}^{X,\lambda}$} The elements of the basis 
are written as $\vec{\Phi}_{j,\gamma,\frakS,n}^\lambda$,
with the following index conventions:
\begin{itemize}
\item 
$j \in \mathbb{J}^\lambda$ is the index related to the background cell,
where for aggregation cells one picks the minimum cell index 
of all aggregated background cells on mesh level $0$ as a representative, 

\item
$\gamma$ is the variable index (e.g. for Navier-Stokes, $1 \leq \gamma \leq D$ corresponds to 
the velocity components and $\gamma = D+1$ corresponds to pressure),

\item
$\frakS \in \{\frakA, \frakB \}$ is the species index and

\item
$1 \leq n \leq N_{k_\gamma}$ is the DG-mode index, where 
$N_k$ is the dimension of the polynomial space up to degree $k$.
\end{itemize}
The row-vector of all XDG basis functions, on mesh level $\lambda$ is written as
\begin{equation}
 \left( 
   \vec{\Phi}_{j,\gamma,\frakS,n}^{X,\lambda}
 \right)_{j,\gamma,\frakS,n}
 =:
 \gvec{\vec{\Phi}}^{X,\lambda} .
\end{equation}
(Here, we skik the specification of all valid combinations of $j,\gamma,\frakS,n$ for sake of compactness.)
In the implementation, this basis 
is constructed from a basis 
$
\underline{\Phi} = \left(
 \Phi_{j,n}
\right)_{ j=1,\ldots,J \atop n=1,\ldots,N_k } 
$
on the space $ \DG{k}(\GridCells) $, with $\supp{\Phi_{j,n}} = K_j$
in the following way:
first, a basis $\underline{\Phi}^\lambda$ of the aggregated space
$ \DG{k}( \aggmesh(\GridCells,A^\lambda) ) $ is created.
Since $\DG{k}( \aggmesh(\GridCells,A^\lambda) ) \leq_\real \DG{k}(\GridCells)$,
one can express this basis in terms of the original basis, i.e. 
\[
  \underline{\Phi}^\lambda = \underline{\Phi} \cdot \gmat{Q}^\lambda,
\]
with a suitable matrix $\gmat{Q}$.
It can be derived from an Ansatz which projects 
a polynomial basis on the bounding box of an aggregation cell onto the 
background cells, for details see ref. \cite{KummerXnse2017}.

The matrix  $\gmat{Q}^\lambda$ obviously is a prolongation operator.
If both, $\underline{\Phi}^\lambda$ and $\underline{\Phi}$ are orthonormal,
the mapping 
\[
  \DG{k}(\GridCells) \ni 
  \underline{\Phi} \cdot \gvec{u} 
    \mapsto
  \underline{\Phi}^\lambda \cdot ((\gmat{Q}^\lambda)^T \gvec{u} ) 
  \in \DG{k}( \aggmesh(\GridCells,A^\lambda) )
\]
is a projector in the $L^2$-sense.

For computational efficiency and in order to save memory, the matrix $\gmat{Q}$ 
should not be stored. Instead, $\underline{\Phi}^{\lambda + 1}$ is expressed in terms 
of $\underline{\Phi}^\lambda$, i.e.
\begin{equation}
  \underline{\Phi}^{\lambda + 1} = \underline{\Phi}^{\lambda} \cdot \gmat{R}^\lambda .
\end{equation}
For orthonormal bases, $\gmat{R}^\lambda$ is the prolongation matrix from the coarse to fine mesh,
while $ ( \gmat{R}^\lambda )^T $ is a restriction, resp. projection matrix from fine to coarse mesh.
For a specific aggregation cell, it can be computed by projection of polynomials 
onto one representative background cell for each part of the aggregation cell.

Second,
the basis of the XDG space is constructed:
the basis elements of $\setV_{\gvec{k}}^{X,\lambda}$ are expressed as 
\begin{equation}
 \vec{\Phi}^{X,\lambda}_{j,\gamma,\frakS,n}  (\vec{x})
  := 
  \vec{e}_{\gamma} \sum_{m=1}^{N_{k_\gamma}} \Phi^\lambda_{j,m}(\vec{x}) \charfunc_{\frakS(t)}(\vec{x}) ~  S_{m n}^\lambda  .
\label{eq:XdgMgBasis}
\end{equation}
Here $\vec{e}_{\gamma}$ is the standard basis vector in $\real^{D_v}$ and $\charfunc_\frakS$ denotes the characteristic function for set $\frakS$.
The matrix $\gmat{S}^\lambda$ provides a re-orthonormalization of the basis functions $\Phi^\lambda_{j,m} \charfunc_{\frakS(t)} $
in cut-cells and can be obtained e.g. through a Cholesky factorization.

Finally, through a combination of multigrid projection matrices $\gmat{R}^\lambda$
and re-orthonormalization matrices $\gmat{S}^\lambda$, one obtains the representation
\begin{equation}
 \gvec{\vec{\Phi}}^{X,\lambda + 1} = \gvec{\vec{\Phi}}^{X,\lambda} \cdot \gmat{R}^{X,\lambda} .
\label{eq:XDGrestrictionMatrix}
\end{equation}
Formally, $\gmat{R}^{X,\lambda}$ is a matrix product of $\gmat{S}^\lambda$ and  $\gmat{R}^\lambda$.
The notation of its exact shape is rather technical and therefore skipped.
Mainly, since e.g. $\gvec{\vec{\Phi}}^{X,\lambda}$ are vector-valued an advanced indexing notation is required,
which is introduced below.

\paragraph{Multi-index mapping}
In order to extract sub-matrices and sub-vectors which correspond to certain 
cells, variables and DG modes, a sophisticated index notation is required.
A single basis element of $\setV_{\gvec{k}}^{X,\lambda}$ can be associated with a multi-index $m(j,\gamma,\frakS,n)$.
One may think of $m(-,-,-,-)$ as a bijection between all valid combinations of 
$j$, $\gamma$, $\frakS$ and $n$ and the set $ I := \{1, \ldots, \dim ( \setV_{\gvec{k}}^{X,\lambda} )  \}$.
We use a notation where the mapping $m(-,-,-,-)$ is employed  to select sub-sets of $I$, e.g.
\[
m(j,-,-,\leq N_{\gvec{k}}) 
:= 
\left\{
   m(j,d_v,\frakS,n) ; \
   \forall d_v, \
   \forall \frakS \text{ present in cell } j, \
   \forall n \leq N_{k_{d_v}}
\right\} .
\]
Such sets will be used to notate sub-matrices or sub-vectors, similar to the typical MATLAB notation.

\paragraph{Agglomeration Algebra}
Given that a basis is established, the generic system (\ref{eq:GenericVarSchemeAgglommSpace}) 
which searches for a solution $U \in  \setV_{\gvec{k}}^{X,\lambda}$
can be transfered to an equivalent matrix formulation
\begin{equation}
 \gmat{M}^\lambda ~ \gvec{U} = \gvec{b}^\lambda
\label{eq:GenericMatrixSystem}
\end{equation}
with 
\begin{equation}
 M^\lambda_{m(j,d_v,\frakS,n) \ m(l,e_v, \frakR, m )} 
 = 
 a( \vec{\Phi}^{X,\lambda}_{l,e_v, \frakR, m}  , \vec{\Phi}^{X,\lambda}_{j,d_v,\frakS,n} )
 \text{  and  }
 b^\lambda_{m(j,d_v,\frakS,n)}  = 
 b( \vec{\Phi}^{X,\lambda}_{l,e_v, \frakR, m} ) .
\end{equation}
Through restriction and prolongation matrices, 
one obtains the relation
\begin{equation}
 \gmat{M}^{\lambda + 1} = ( \gmat{R}^{X,\lambda} )^T ~ \gmat{M}^{\lambda} ~ \gmat{R}^{X,\lambda}
 \text{ and }
 \gvec{b}^{\lambda + 1} = ( \gmat{R}^{X,\lambda} )^T ~ \gvec{R}^{\lambda} ,
\label{eq:XDGsysRestriction}
\end{equation}
as usual in multigrid methods.
% ===============================================================================
\subsection{Preconditioners and Solvers}
\label{sec:Precond}
% ===============================================================================
Within this section the discussion is focused on a single mesh level $\lambda$, hence the 
level index is dropped.
On this grid level, let $L$ be the vector space dimension, i.e.
$\gmat{M} = \gmat{M}^\lambda \in \real^{L \times L}$.
The (approximate) solution and right-hand-side (RHS) vector 
of the system (\ref{eq:GenericMatrixSystem}) are denoted as 
$\gvec{x}, \gvec{b} \in \real^L$, respectively.

In this section, a series of algorithms is introduced.
\begin{itemize}
\item 
A p-multigrid algorithm (Algorithm \ref{alg:pMultigrid}), that operates on a single mesh level,
which can be used as a pre-conditioner for a standard GMRES solver.

\item 
The additive Schwarz algorithm (Algorithm \ref{alg:Swz}), which uses p-multigrid as a block solver.
It contains a minor modification to its original form, which we found helpful in 
decreasing the number of iterations and is therefor also presented here.

\item 
Finally, the orthonormalization multigrid algorithm (Algorithm \ref{alg:OrthoMG})
which employs the additive Schwarz as a smoother and a residual minimization (Algorithm \ref{alg:ResMinimi})
to ensure non-increasing residual norms.
\end{itemize}

\paragraph{A p-multigrid pre-conditioner on a single mesh level}
For DG or XDG methods, even without aggregation meshes, a p-multigrid 
seems to be a reasonable idea: At first, a sub-matrix which corresponds to 
low order DG modes is extracted. Since the number of degrees-of-freedom for this is 
typically low, a sparse direct solver can be used.
The degrees of freedom which correspond to higher order DG modes are then solved locally 
in each cell. Since Algorithm \ref{alg:pMultigrid} will also be used as a block solver for 
the additive Schwarz method (Algorithm \ref{alg:Swz}), an optional 
index-set $\mathbb{I} \subset \mathbb{J}^\lambda$ which restricts the solution to a part of the mesh
can be provided.
\begin{myalgorithm}[p-multigrid pre-conditioner]
\label{alg:pMultigrid}
$\gvec{x} = \text{PMG}(\gvec{b},\gvec{k}_{\text{lo}},\mathbb{I})$  is computed as follows: 
\begin{algorithmic}
\Require A right-hand-side $\gvec{b}$,
         a vector of DG polynomial degrees $\gvec{k}_{\text{lo}}$ 
         which separates low from high order modes, for each variable.
         Optionally, an index-set $\mathbb{I} \subset \mathbb{J}^\lambda$ of 
         sub-cells, denoting the block to solve.
\Ensure  An approximate solution $\gvec{x}$
\State
$\gvec{x} := 0$ 
\Comment{initialize approximate solution}
\State $I_{\text{lo}} := m(\mathbb{I},-,-,\leq N_{\gvec{k}_{\text{lo}}})$
\Comment{Indices for low-order modes in cells $\mathbb{I}$}
\State 
$\gmat{M}_{\text{lo}} := \gmat{M}_{ I_{\text{lo}} , I_{\text{lo}} }$
\Comment{extract low-order system}
\State 
$\gvec{b}_{\text{lo}} := \gvec{b}_{ I_{\text{lo}} }$
\Comment{extract low-order RHS}
\State Solve: $\gvec{x}_{\text{lo}} := \gmat{M}_{\text{lo}} \setminus \gvec{b}_{\text{lo}}$
\Comment{usually using a sparse direct solver}
\State
$\gvec{x}_{ I_{\text{lo}} } := \gvec{x}_{ I_{\text{lo}} } + \gvec{x}_{\text{lo}}$
\Comment{accumulate low-order solution}
\State
$\gvec{r} := \gvec{b} - \gmat{M} \gvec{x} $
\ForAll {$j \in \mathbb{I}$}
\Comment{loop over cells...}
  \State $I_{\text{hi},j} := m(j,-,-,> N_{\gvec{k}_{\text{lo}}})$
  \Comment{Indices for high-order modes in cell $j$}
  \State $\gmat{M}_{\text{hi},j} := \gmat{M}_{I_{\text{hi},j} \ I_{\text{hi},j}}$
  \Comment{extract high-order system in cell $j$}
  \State $\gvec{b}_{\text{hi},j} := \gmat{b}_{ I_{\text{hi},j} }$
  \Comment{extract high-order RHS in cell $j$}
  \State Solve: $\gvec{x}_{\text{hi},j} := \gmat{M}_{\text{hi},j} \setminus \gvec{b}_{\text{hi},j}$
  \Comment{using a dense direct solver}
  \State $\gvec{x}_{ I_{\text{hi},j} } := \gvec{x}_{ I_{\text{hi},j} } + \gvec{x}_{\text{hi},j}$
  \Comment{accumulate local high-order solution}
\EndFor
\end{algorithmic}
\end{myalgorithm}
Since a preconditioner typically has to be applied multiple time,
it is essential for performance to use a sparse direct solver for the low order system 
which is able to store the factorization and apply multiple right-hand sides.
Numerical tests hint that usually a single-precision solvers is sufficient 
as long as the residuals are computed in double precision.
The PMG aglorithm presented above can directly be used as a preconditioner, e.g. for GMRES.
However, since only two multigrid levels are used, its application typically 
is only reasonable up to medium-sized systems.

\paragraph{An additive Schwarz method}
In the scope of this work, $\text{PMG}(\ldots)$ is also used as a block solver for an additive Schwarz method.
For a general discussion on Schwarz methods, we refer to the review article of Gander \cite{Gander2008}.
In our implementation the blocking is determined on the level of cells, i.e. on the basis of the 
mesh $\GridCells^\lambda$. The METIS \cite{Karypis1998} software library is used to partition $\GridCells^\lambda$;
the number of partitions is determined so that a single partition contains roughly 10,000 degrees-of-freedom.
After the partitioning by METIS, each partition is enlarged by its neighbor cells to generate the overlap
layer that is typically used with Schwarz methods.
\begin{myalgorithm}[Additive Schwarz] $\gvec{x} = \text{Swz}(\gvec{b})$ is computed as follows: 
\label{alg:Swz}
\begin{algorithmic}
\Require A right-hand-side $\gvec{b}$,
         a vector of DG polynomial degrees $\gvec{k}_{\text{lo}}$ 
         which separates low from high order modes, for each variable.
         Furthermore, a pre-computed partitioning of cells, i.e. sets $\mathbb{I}_i \subset \mathbb{J}_\lambda$
         so that $\cup_i \mathbb{I}_i = \mathbb{J}^\lambda$.
\Ensure  An approximate solution $\gvec{x}$
\State
$\gvec{x} := 0$ 
\State
$\gvec{\alpha} := 0$
\Comment{initialize approximate solution}
\State
$\gvec{r} := \gvec{b} - \gmat{M} \gvec{x} $
\ForAll {$\mathbb{I}_i$} 
\Comment{Loop over Schwarz blocks...}
  \State $\gvec{x} := \gvec{x} + \text{PMG}(\gvec{r},\gvec{k}_{\text{lo}},\mathbb{I}_i)$
  \Comment{solve block no. $i$}
  \State $\gvec{\alpha} := \gvec{\alpha} + \charfunc_{m(\mathbb{I}_i,-,-,-)}$
\EndFor
\State $\gvec{x} := \gvec{x} ./ \gvec{\alpha}$
\Comment{Apply scaling in overlap regions}
\end{algorithmic}
\end{myalgorithm}
Within the scope of this work, the additive Schwarz solver is used as a smoother for 
the multigrid method. 
As such, we introduced a minor modification to the original formulation:
in the last line of Algorithm \ref{alg:Swz}, the solution is divided by the number of blocks which 
contain a specific cell. This damping of the approximate solution in the overlapping regions has 
shown to improve the number of iterations when $\text{Swz}(-)$ is used as a smoother for 
the multigrid method (Algorithm \ref{alg:OrthoMG}).

\paragraph{A multigrid algorithm based on orthonormalization}
Experiments have shown that the classic multigrid algorithm,
as could be found e.g. in the textbook of Saad \cite{SaadBook2003} often diverges for DG and XDG methods.
To overcome this issue, we employ a residual minimization approach:
Let $\gmat{Z} = (\gvec{z}_1,\ldots,\gvec{z}_{l}) \in \real^{L \times l}$ be 
a family of approximate solutions, e.g. from some Additive Schwarz or from a coarse grid solver.
Furthermore let be 
$\gmat{W} = (\gvec{w}_1,\ldots,\gvec{w}_{l}) = \gmat{M} ~ \gmat{Z} \in \real^{L \times l}$  
be the image of $\gmat{Z}$ under the application of $\gmat{M}$.
For the residual-optimized solution, one employs the Ansatz
\begin{equation}
  \gvec{x} = \gvec{x}_0 + \sum_i \alpha_i \gvec{z}_i ,
\end{equation}
in order to  minimize the 2-norm of the residual
\begin{eqnarray}
 \gvec{r} & := & \gvec{b} - \gmat{M} ~ \gvec{x} \\
          &  = &\gvec{r}_0 - \gmat{M} ~ \gmat{Z} ~ \gvec{\alpha} \\
          &  = & \gvec{r}_0 - \gmat{W} ~ \gvec{\alpha} .
\end{eqnarray}
As widely known, 
the coefficients $\gvec{\alpha} $ at which $| \gvec{r} |_2$ becomes minimal are determined by 
the system 
\begin{equation}
  \left( \gmat{W}^T \cdot \gmat{W} \right) \gvec{\alpha} = \gmat{W}^T \cdot \gvec{r}_0 .
\label{eq:ResMinimi}
\end{equation}
Furthermore, if columns of $\gmat{W}$ are orthonormal, the left-hand-side of equation (\ref{eq:ResMinimi})
becomes identity. This can be achieved through a Gram-Schmidt algorithm, which is used in the following
residual minimization algorithm which is the foundation of our modified multigrid method:
\begin{myalgorithm}[Residual minimization]
\label{alg:ResMinimi}
$(\gvec{x}, \gvec{r}) = \text{RM}(\gvec{x}_0, \gvec{z}, \gmat{W}, \gmat{Z})$  is computed as follows: 
\begin{algorithmic}
\Require 
A solution-guess $\gvec{z} \in \real^L$, e.g. from some pre-conditioner.
An orthonormal system 
$\gmat{W} = (\gvec{w}_1,\ldots,\gvec{w}_{l}) \in \real^{L \times l}$
and its application onto the system matrix 
$\gmat{Z}$, i.e. $\gmat{W} = \gmat{M}  \gmat{Z}$.
Note that $l = 0$ is allowed, i.e. $\gmat{W}$ and $\gmat{Z}$ can be empty.
\Ensure  An approximate solution $\gvec{x}$ whose residual $\gvec{r}$ is not greater in the 
$l_2$-norm than the residual of the initial value $\gvec{x}_0$.
Furthermore, updated $\gmat{W}, \gmat{Z} \in \real^{L \times l+1}$ 
where the columns of $\gmat{W}$ are orthonormal and $\gmat{W} = \gmat{M}  \gmat{Z}$.
\State{$\gvec{r}_0 = \gvec{b} - \gmat{M} \gvec{x}_0 $ }
\State{$\gvec{w}_{l+1} = \gmat{M} \gvec{z}$}
\ForAll{columns $\gvec{w}_{i}$ of $\gmat{W}$}
\Comment{Gram-Schmidt loop...}
  \State{$\beta := \gvec{w}_i \cdot \gvec{w}_{l+1}$}
  \State{$\gvec{w}_{l+1} := \gvec{w}_{l+1} - \beta \gvec{w}_i$}
  \State{$\gvec{z} := \gvec{z} - \beta \gvec{Z}_i$}
\EndFor
\State{$\gamma := 1 / | \gvec{w}_{l+1}  |_2 $, 
	   $\gvec{w}_{l+1} := \gamma \gvec{w}_{l+1}$, 
	   $\gvec{z}_{l+1} := \gamma \gvec{z}_{l+1}$}
\State{$\gmat{W} := ( \gmat{W} , \gvec{w}_{l+1} ) $, 
       $\gmat{Z} := ( \gmat{Z} , \gvec{z} )$      }
\Comment{Store vectors $\gvec{w}_{l+1}$ and $\gvec{z}$}
\State{$\alpha_{l+1} := \gvec{w}_{l+1} \cdot \gvec{r}_0$, $\gvec{\alpha} := (\gvec{\alpha}, \alpha_{l+1})  $ }
\State{$\gvec{x} := \gvec{x}_0 + \gmat{Z} \cdot \gvec{\alpha} $}
\Comment{optimized solution}
\State{$\gvec{r} := \gvec{r}_0 - \gmat{W} \cdot \gvec{\alpha} $}
\Comment{residual for optimized solution}
\end{algorithmic}
\end{myalgorithm}

In order to define a multigrid algorithm, one employs the 
restriction matrix $\gmat{R}^{X,\lambda} \in \real^{L \times L_c}$
introduced in Eq. (\ref{eq:XDGrestrictionMatrix}).
from which one infers the restricted-system matrix $\gmat{M}^{\lambda+1}$, cf. Eq. (\ref{eq:XDGsysRestriction}).
Here, $L_c$ denotes the vector-space dimension of the coarser space  $\setV_{\gvec{k}}^{X,\lambda+1}$
at mesh level $\lambda + 1$.
As noted, the classical multigrid algorithm very often does not converge for DG or XDG methods,
at least if Block-Jacobi or Schwarz algorithms are used a smoothers.
Therefore, we employ $\text{RM}(\ldots)$ after each pre- and post-smoother step and after  
the coarse-grid correction in order to ensure a non-increasing residual norm.
	
\begin{myalgorithm}[Orthonormalization multigrid]
\label{alg:OrthoMG}
$(\gvec{x}, \gvec{r}) = \text{MG}(\gvec{x}_0, \gvec{b}, \gmat{M})$  is computed as follows: 
\begin{algorithmic}
\Require A right-hand-side $\gvec{b}$ and an initial solution guess $\gvec{x}_0$ (can be zero).
\Ensure An approximate solution $\gvec{x}$ whose residual norm is 
       less than or equal to the residual norm of the initial guess $\gvec{x}_0$. 
\State $\gmat{Z} := (), \gmat{W} := ()$
\Comment{Initialize as empty}
\State $\gvec{r}_0 := \gvec{b} - \gmat{M} \gvec{x}_0$, $\gvec{r} := \gvec{r}_0$, $\gvec{x} := \gvec{x}_0$
\Comment{residual of interstitial solution}
\While{$|\gvec{r}|_2 > \varepsilon$}
  \State{ $\gvec{z} := \text{Swz}(\gvec{b})$}
  \Comment{pre-smoother}
  \State{$(\gvec{x}, \gvec{r}) = \text{RM}(\gvec{x}_0, \gvec{z}, \gmat{W}, \gmat{Z})$}
  \Comment{minimize residual of pre-smoother}
  \State{$\gvec{r}_c := (\gmat{R}^{X,\lambda})^T \gvec{r}$}
  \Comment{restrict residual}
  \State{$(\gvec{z}_c, \gvec{r}_c) = \text{MG}(0,\gvec{r}_c, \gmat{M}^{\lambda + 1})$}
  \Comment{call multigrid on coarser level}
  \State{$\gvec{z} := \gmat{R}^{X,\lambda} \gvec{z}_c $}
  \Comment{prolongate coarse-grid correction}
  \State{$(\gvec{x}, \gvec{r}) = \text{RM}(\gvec{x}_0, \gvec{z}, \gmat{W}, \gmat{Z})$}
  \Comment{minimize residual}
  \State{ $\gvec{z} := \text{Swz}(\gvec{b})$}
  \Comment{pre-smoother}
  \State{$(\gvec{x}, \gvec{r}) = \text{RM}(\gvec{x}_0, \gvec{z}, \gmat{W}, \gmat{Z})$}
  \Comment{minimize residual of pre-smoother}
\EndWhile
\end{algorithmic}
\end{myalgorithm}

% ===============================================================================
% ===============================================================================
\section{Solver Performance studies}
\label{chap:SolverPerformance}
% ===============================================================================
% ===============================================================================
The performance of solvers introduced in section \ref{sec:Precond} is investigated
for the XDG Poisson problem (\ref{eq:poisson-jump-problem-def}) 
on the domain $\Omega = (-1,1)^3$
with homogeneous Dirichlet boundary conditions ($\GammaDiri = \partial \Omega$, $g_\text{Diri} = 0$), right hand side $f=1$
and a large ratio of diffusion coefficients ($\mu_\frakA = 1$, $\mu_\frakB = 1000$)
for a Level-Set $\varphi = x^2 + y^2 + z^3 - (7/10)^2$.
Equidistant, Cartesian meshes were employed with a resolution of 
$2^3$, $4^3$, $8^3$, etc., cells and  polynomial degrees $k$ of $2$, $3$ and $5$.
The agglomeration threshold  $\alpha$ (cf. section \ref{sec:Agglom})
is set 0.1 for $k=2,3$ and 0.3 for $k=5$.
Details on degrees-of-freedom can be found in Table \ref{tab:DOFs}.
All solvers were configured to terminate if the 2-norm of residual is below $10^{-10}$.

In order to assess runtime-measurements, one first has to ensure that the 
subroutines of which the iterative solvers are composed of.
The most complex of those a 
the sparse matrix-vector and matrix-matrix product. A (non-systematic) 
comparison against MATLAB shows that the implementations in BoSSS are
quite decent, as can bee seen in Table \ref{tab:MatrixOpBench}.
The test matrix for this benchmark 
is the XDG system used in this section, for $6^3$ cells and a polynomial order 
$k = 5$, yielding a matrix with 24,192 non-zero rows. It should be noted that 
the BoSSS implementation exploits the block-structure of the XDG matrix, 
thus BLAS (Intel MKL 11.0) subroutines can be used for the inner dense 
matrix-matrix an matrix-vector products. Since BoSSS is MPI parallel, only 
the sequential version of the MKL is used in order to avoid that all MPI 
ranks try to use all cores of some compute node. 

\begin{figure}[!h]
\begin{center}
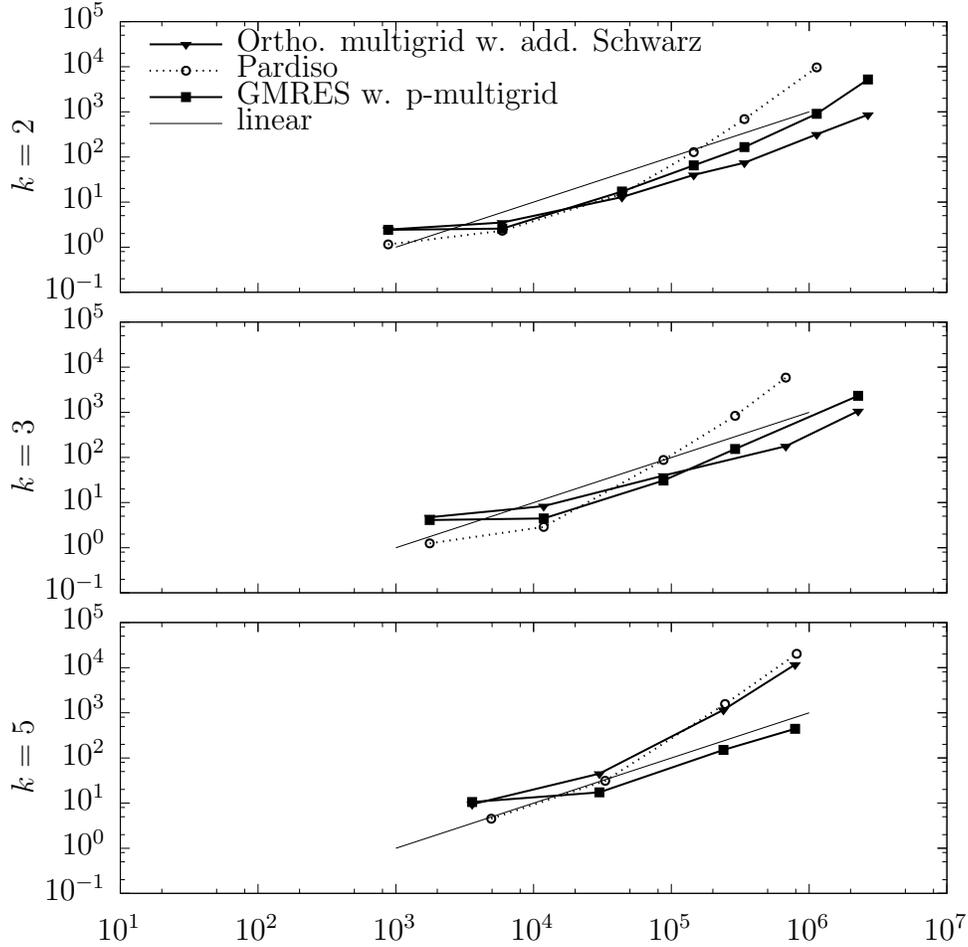
\end{center}
\caption{
Overall solver runtime (in seconds) versus. degrees-of-freedom, for different polynomial degrees $k$,
for the XDG Poisson problem with 1:1000 diffusion coefficient ratio.
The solver are: PARDISO, orthonormalization multigrid (alg. \ref{alg:OrthoMG}), GMRES with p-multigrid (alg. \ref{alg:pMultigrid}).
}
\label{fig:XdgRuntimes}
\end{figure}

Results of run-time measurements can be seen in Figure \ref{fig:XdgRuntimes}.
A performance base-line is given by the runtime of 
the (direct sparse) PARDISO solver \cite{Schenk2002,Schenk2004,Schenk2006}.
As expected, direct sparse solvers are superior for comparatively small problems
up 100,000 degrees-of-freedom. In addition to the presented results, 
we also tested the MUMPS package \cite{AmestoyMUMPS2001, AmestoyMUMPS2006}, showing only 
minor differences in runtime behavior, which are therefore skipped.
Despite having a symmetric, positive definite system the direct solver was configured to 
assume general matrices, since we are finally interested in non-symmetric systems such as the 
Navier-Stokes problem (\ref{eq:TwoPhaseNSE}).
Beyond 100,000 degrees-of-freedom, iterative solvers out-perform sparse direct systems
and show a linear run-time behavior, with one exception.

\begin{figure}[!h]
\begin{center}
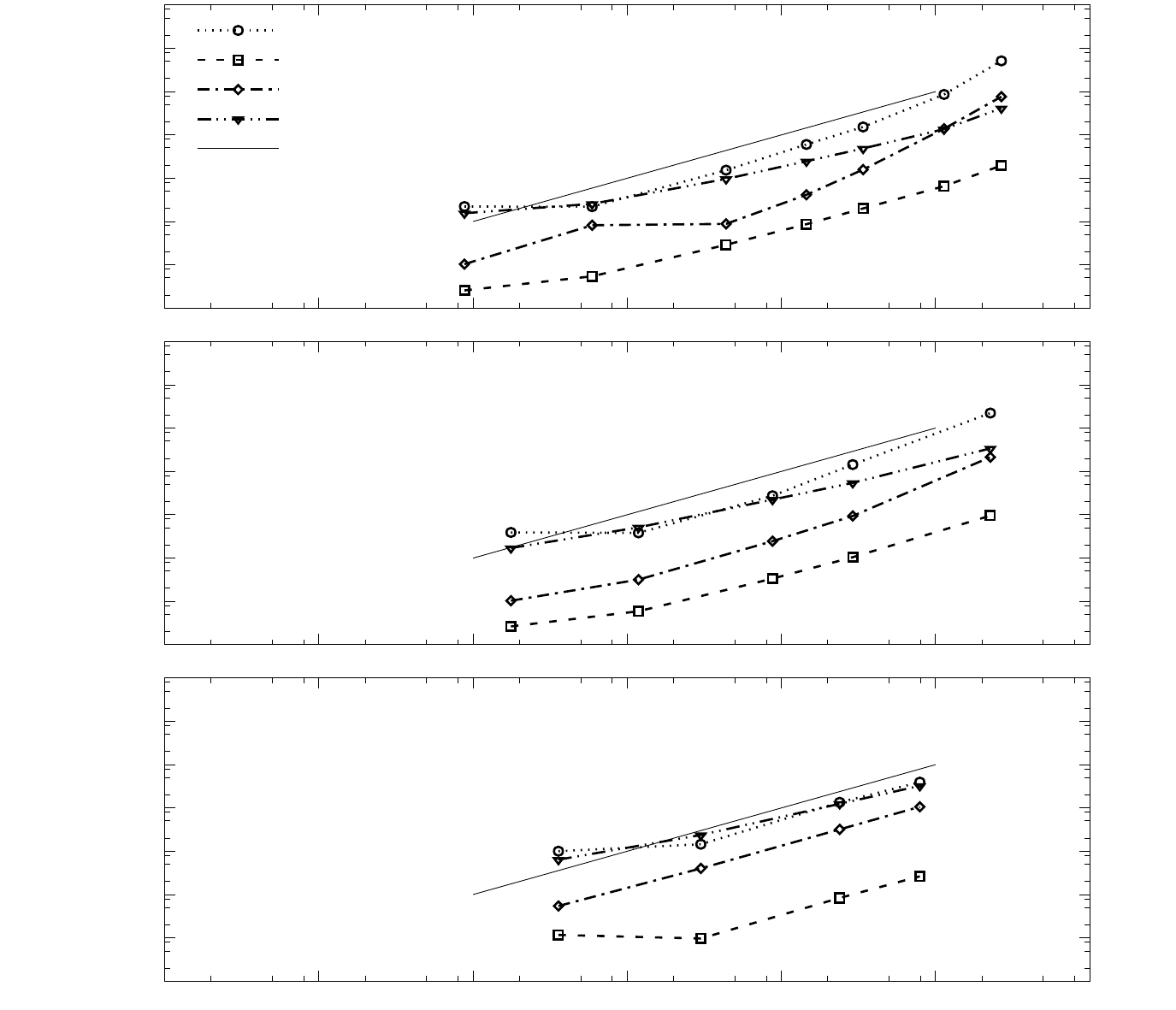
\end{center}
\caption{
Runtime (in seconds) of individual solver phases of the GMRES with p-multigrid (alg. \ref{alg:pMultigrid})
 versus degrees-of-freedom for different polynomial degrees $k$,
for the XDG Poisson problem with 1:1000 diffusion coefficient ratio.
Note that the matrix assembly is not counted as part of the solution process and is only given for purpose additional perspective.
}
\label{fig:GMRES_PMG}
\end{figure}

\begin{figure}[!h]
\begin{center}
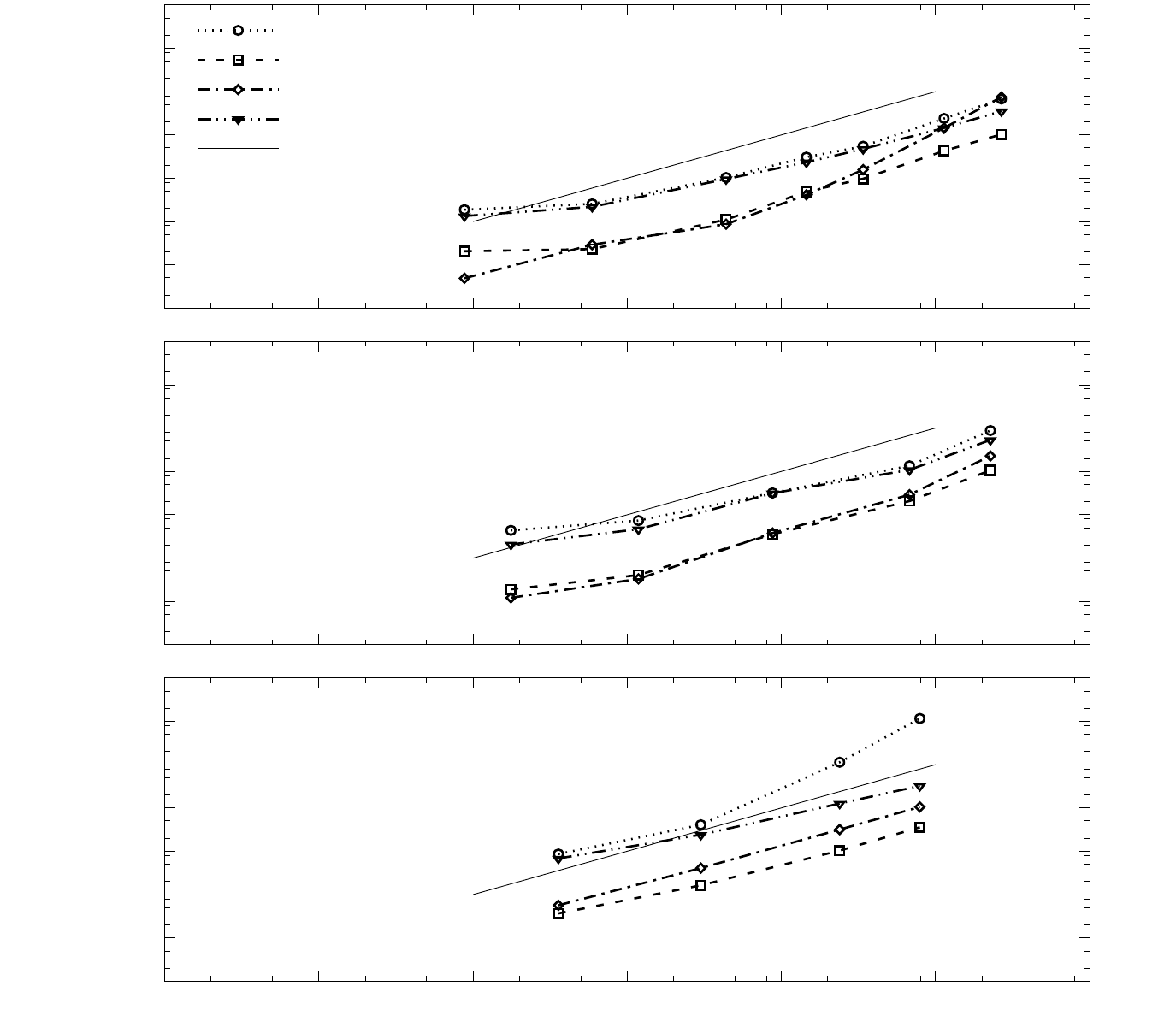
\end{center}
\caption{
Runtime (in seconds) of individual solver phases of the orthonormalization multigrid (alg. \ref{alg:OrthoMG})
 versus degrees-of-freedom for different polynomial degrees $k$,
for the XDG Poisson problem with 1:1000 diffusion coefficient ratio.
Note that the matrix assembly is not counted as part of the solution process and is only given for purpose additional perspective.
}
\label{fig:Kcycle_Schwarz}
\end{figure}

Runtime-measurements also include setup-costs for the iterative solvers.
Figures \ref{fig:Kcycle_Schwarz} and \ref{fig:GMRES_PMG} outline 
runtime details for the orthonormalization multigrid (Algorithm \ref{alg:OrthoMG})
as well as GMRES with p-multigrid preconditioner (Algorithm \ref{alg:pMultigrid}), respectively.
As individual solver phases, we identify 
(i)
the setup of the XDG aggregation grid basis, including the re-orthonormalization 
in cut-cells (matrix $\gmat{S}^\lambda$, cf. Eq. \ref{eq:XdgMgBasis}),
(ii) 
the sparse matrix-matrix products required to setup the sequence of 
linear systems (cf. Eq. \ref{eq:XDGsysRestriction}) and
(iii)
the iterative procedure itself.

A clearly negative result here is the worse-than-linear runtime behavior 
for the orthonormalization multigrid for $k=5$, cf. Figure \ref{fig:XdgRuntimes}.
While this is caused by the solver iterations (cf. Figure \ref{fig:Kcycle_Schwarz}),
the number of iterations  is fairly constant for the number of degrees-of-freedom.
In order to achieve this, the separation degree $k_{\text{lo}}$ (cf. Algorithm \ref{alg:pMultigrid})
had to be increased from 1 for the cases $k=2,3$ to 3 for the case $k=5$.
This hints (an is also verified by measurements)  that the time for a single iteration 
is dominated by calls to the sparse linear solver for the low-order Schwarz blocks.
We want to point out that such a behavior could not be observed using a standard DG 
method - there, $k_{\text{lo}} = 1$ seems to be sufficient, whith only a minor dependence 
of the number of iterations on $k$.
Ongoing works investigate an adaptation of $k_{\text{lo}}$ in cut-cells, in order to 
regain the runtime behavior for DG problems for XDG.

\begin{table}
	\begin{center}
		\begin{tabular}{l||c|c|c|c|c|c|c|c}
			&\multicolumn{ 8}{|c}{grid resolution}\\
			%\cline{2-9}
			& $2^3$ &	$4^3$ &	$8^3$	 & $16^3$ &	$24^3$ &	$32^3$ &	$48^3$ & $64^3$ \\
			%\cline{2-9}
			%& 8 &	64 &	512	 & 4096 &	13824 &	32768 &	110592 & 262144 \\
			\hline \hline
			$p=2$ & 160	& 880  & 5920  &	43840 &	145200	 &339040 &	1134560	 & 2671600 \\
			%\hline
			$p=3$ & 320	& 1760	& 11840	 &87680 &	290400 &	678080	 & 2269120 &	\\
			%\hline
			$p=5$ & 896	& 4928	& 33152 &	245504 &	813120  &  &  &  \\
		\end{tabular}
	\end{center}
\caption{degrees of freedom in dependency of grid resolution and polynomial-degree $k$ for the subsequent Performance plots shown in figure \ref{fig:XdgRuntimes}, \ref{fig:XdgIter}, \ref{fig:GMRES_PMG} and \ref{fig:Kcycle_Schwarz} for the problem described in section \ref{chap:SolverPerformance}}
\label{tab:DOFs}
\end{table}

\begin{table}
\begin{center}
\begin{tabular}{l||c|c|c}
                            & \multicolumn{3}{c}{matrix-matrix}    \\ 
                            & BoSSS  & BoSSS   & Matlab              \\ 
                            & 1 MPI core & 4 MPI cores & OpenMP   \\ 
\hline \hline
i7 6700HQ (4 cores)         & 3.88  &  1.293   & 5.19                \\
Xeon E5-2630L v4 (20 cores) & 5.86  &  1.687   & 4.42                \\            
\multicolumn{4}{c}{}                                                 \\
                            & \multicolumn{3}{c}{matrix-vector}      \\ 
\hline \hline
i7 6700HQ (4 cores)         & $14.6 \cdot 10^{-3}$ &  $8.99  \cdot 10^{-3}$  &  $11.8 \cdot 10^{-3}$ \\ 
Xeon E5-2630L v4 (20 cores) & $29.2 \cdot 10^{-3}$ &  $9.14  \cdot 10^{-3}$  &  $7.53 \cdot 10^{-3}$ \\  
\end{tabular}
\end{center}
\caption{
Runtime (in seconds) of a sparse matrix-matrix and matrix-vector multiplication
in BoSSS vs. MATLAB on two different computer systems for a quadratic matrix with 24,192 non-zero rows.
This verifies that these building blocks provide decent performance,
given that the BoSSS implementation is sequential (within each MPI core) and MATLAB is OpenMP parallel.
An advantage of the BoSSS implementation is that it exploits the block-structure of the XDG matrix, which is not known to  MATLAB.
} 
\label{tab:MatrixOpBench}
\end{table}

\begin{figure}[!h]
\begin{center}
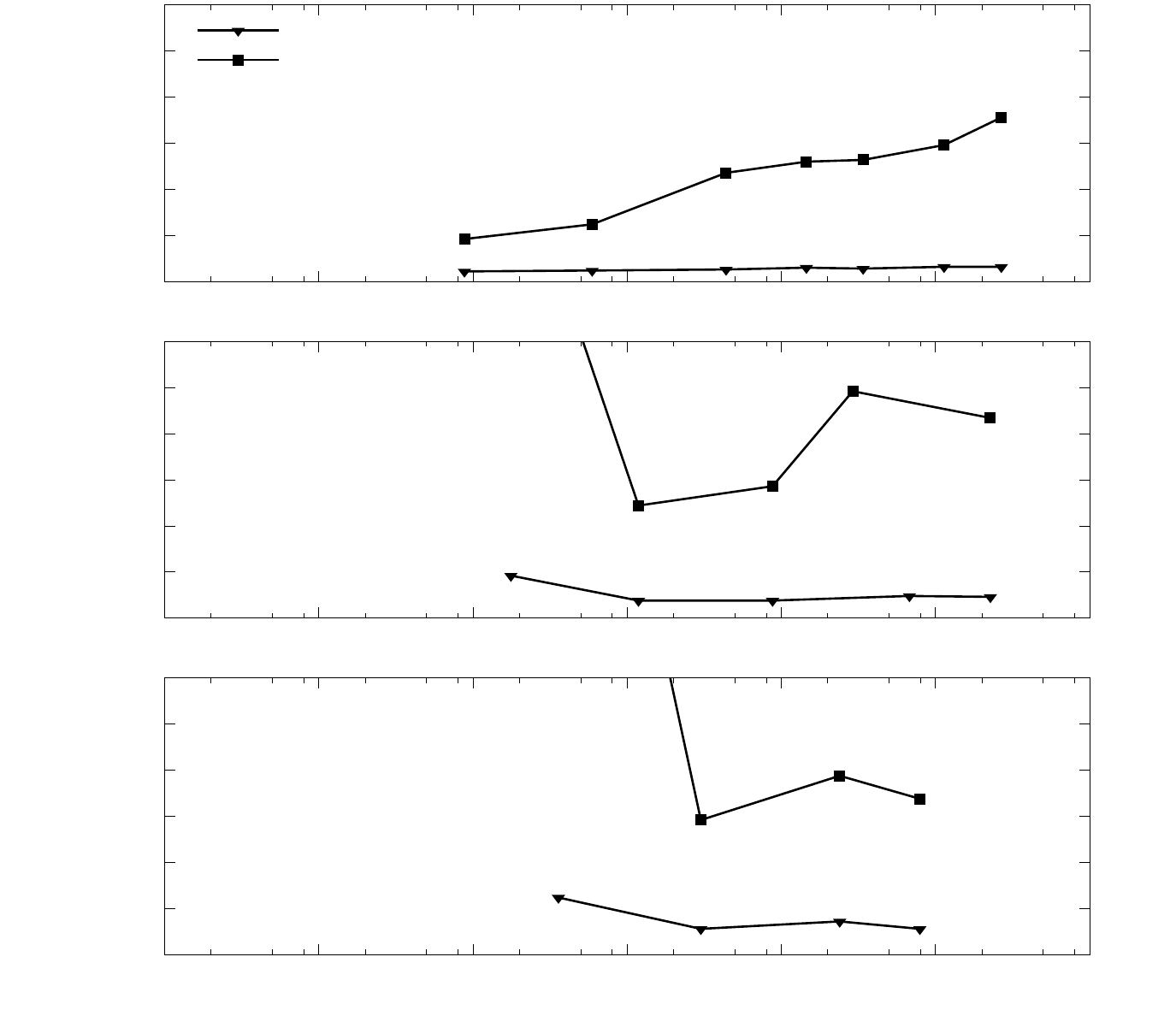
\end{center}
\caption{
Number of iterations versus degrees-of-freedom for different polynomial degrees $k$ 
for the XDG Poisson problem with 1:1000 diffusion coefficient ratio.
The number of iterations for the orthonormalization multigrid is nearly independent 
from the number of degrees-of-freedom.
GMRES with a single level of p-multigrid shows a slight increase of iterations with 
outliers at the lower end, where direct solvers are superior.
Despite the higher number of iterations, GMRES might still be faster in wall-clock time
cf. Figure \ref{fig:XdgRuntimes}.
}
\label{fig:XdgIter}
\end{figure}

% ===============================================================================
% ===============================================================================
\section{Conclusion and Outlook}
\label{sec:Conclusions}
% ===============================================================================
% ===============================================================================

In this work, the XDG framework available in BoSSS for the discretization 
of PDEs with discontinuous coefficients was presented.
and
a unified notation for the agglomeration of small cut-cells and aggregation multigrid
was developed.
The practical advantages of XDG for single phase flows is
the elimination of the meshing process, which is especially useful for time-evolving geometries.
Regarding two-, resp. multiphase flows, the XDG method preserves the $h^{k+1}$-convergence rate 
despite the presence of discontinuous material coefficients, cf. \cite{KummerXnse2017,MuellerEtAl2017}.

Furthermore, in this work we presented a combined p- and h-multigrid method and assessed its performance.
Althogh the BoSSS code is fully MPI parallel, up to this point, parallel scaling 
has not been investigated \emph{systematically}
and only  set of various results such as Table \ref{tab:MatrixOpBench} exist. 
Extensive tests for parallel efficiency of the presented solvers and other subsystems of BoSSS 
are a matter of ongoing work and will be published in the BoSSS handbook, 
available at the online code repository at Github.

So far, only the solution of linear systems was addressed.
Newton-methods for the nonlinear Navier-Stokes problem (\ref{eq:TwoPhaseNSE}) coupled together with the 
Level-Set evolution (\ref{eq:LevelSetEq}) are subject of ongoing works and will 
be the issue of upcoming publications.

% ===============================================================================
% ===============================================================================
\section*{Acknowledgment} 
% ===============================================================================
% ===============================================================================
The work of F. Kummer and M. Smuda is supported by 
the German Science Foundation (DFG) through 
the Collaborative Research Centre (SFB) 1194
``Interaction between Transport and Wetting Processes''.



%% References


\bibliographystyle{elsarticle-num}
\bibliography{OpenSoftwarePDE}

\end{document}